\documentclass[5p,twocolumn]{elsarticle} 

\journal{Automatica}
\biboptions{numbers,sort&compress}

\usepackage{amsmath,amssymb,amsthm}
\usepackage{algorithm,algorithmic,setspace}
\usepackage[hidelinks]{hyperref}
\usepackage{bbm}

\usepackage[shortlabels]{enumitem}

\usepackage{mathtools}
\newcommand{\eqlabel}[1]{\stackrel{\mathclap{\mathrm{#1}}}{=}}

\newcommand{\geqlabel}[1]{\stackrel{\mathclap{\mathrm{#1}}}{\geq}}

\usepackage{tikz}

\newcommand{\algspc}{1.2}       
\allowdisplaybreaks

\newtheorem{theorem}{Theorem}
\newtheorem{proposition}{Proposition}
\newtheorem{lemma}{Lemma}
\newtheorem{assumption}{Assumption}

\DeclareMathOperator{\blkdiag}{blkdiag}
\DeclareMathOperator*{\argmin}{argmin}

\newcommand{\T}{^\top}
\newcommand{\norm}[1]{\|#1\|}
\newcommand{\map}[3]{#1: #2 \rightarrow #3}
\newcommand{\bs}{\boldsymbol}
\newcommand{\until}[1]{\{1,\dots,#1\}}
\newcommand{\toN}{1,\dots,\Nag}

\newcommand{\bone}{\mathbbm{1}}
\newcommand{\bN}{\mathbb{N}}
\newcommand{\bR}{\mathbb{R}}

\newcommand{\cD}{\mathcal{D}}
\newcommand{\cE}{\mathcal{E}}
\newcommand{\cG}{\mathcal{G}}
\newcommand{\cP}{\mathcal{P}}
\newcommand{\cV}{\mathcal{V}}
\newcommand{\nbrs}{\mathcal{N}}

\newcommand{\dual}{\varphi}

\newcommand{\kron}{\otimes}

\newcommand{\cW}{\mathcal{W}}
\newcommand{\tW}{\tilde{W}}
\newcommand{\ones}{\mathcal{O}}

\newcommand{\cC}{\mathcal{C}}

\newcommand{\iter}{k}
\newcommand{\Nag}{N}

\newcommand{\cK}{\mathcal{K}}
\newcommand{\seq}[1]{{\{#1\}}_{\iter \geq 0}}
\newcommand{\subseq}[1]{\{#1\}_{\iter \in \cK}}
\newcommand{\xlim}{\bs{\tilde{x}}}
\newcommand{\lavglim}{\tilde{\lambda}}

\newcommand{\xik}{x_{i,\iter}}
\newcommand{\xikp}{x_{i,\iter +1}}
\newcommand{\xis}{x_i^\star}

\newcommand{\lik}{\lambda_{i,\iter}}
\newcommand{\ljk}{\lambda_{j,\iter}}
\newcommand{\likp}{\lambda_{i,\iter +1}}
\newcommand{\ls}{\lambda^\star}
\newcommand{\dik}{d_{i,\iter}}
\newcommand{\djk}{d_{j,\iter}}
\newcommand{\dikp}{d_{i,\iter +1}}
\newcommand{\deltaik}{\delta_{i,\iter}}
\newcommand{\ellik}{\ell_{i,\iter}}
\newcommand{\deltajk}{\delta_{j,\iter}}
\newcommand{\elljk}{\ell_{j,\iter}}

\newcommand{\xx}{\bs{x}}
\renewcommand{\ll}{\bs{\lambda}}
\newcommand{\dd}{\bs{d}}
\newcommand{\zz}{\bs{z}}
\newcommand{\ee}{\bs{e}}

\newcommand{\xk}{\xx_{\iter}}
\newcommand{\xkp}{\xx_{\iter +1}}
\newcommand{\xs}{\xx^\star}
\newcommand{\lk}{\ll_{\iter}}
\newcommand{\lkp}{\ll_{\iter +1}}
\newcommand{\lls}{\ll^\star}
\newcommand{\dk}{\dd_{\iter}}
\newcommand{\dkp}{\dd_{\iter +1}}
\newcommand{\zk}{\zz_{\iter}}
\newcommand{\zkp}{\zz_{\iter +1}}
\newcommand{\deltak}{\bs{\delta}_{\iter}}
\newcommand{\ellk}{\bs{\ell}_{\iter}}

\newcommand{\dkavg}{\bar{d}_{\iter}}
\newcommand{\dkpavg}{\bar{d}_{\iter +1}}
\newcommand{\lkavg}{\bar{\lambda}_{\iter}}
\newcommand{\lkpavg}{\bar{\lambda}_{\iter +1}}
\newcommand{\ddkavg}{\bs{\bar{d}}_{\iter}}
\newcommand{\ddkpavg}{\bs{\bar{d}}_{\iter +1}}
\newcommand{\llkavg}{\bs{\bar{\lambda}}_{\iter}}
\newcommand{\llkpavg}{\bs{\bar{\lambda}}_{\iter +1}}

\newcommand{\eed}{\ee^{\dd}}
\newcommand{\eel}{\ee^{\ll}}
\newcommand{\edk}{\eed_{\iter}}
\newcommand{\edkp}{\eed_{\iter +1}}
\newcommand{\elk}{\eel_{\iter}}
\newcommand{\elkp}{\eel_{\iter +1}}
\newcommand{\ek}{\ee_{\iter}}
\newcommand{\ekp}{\ee_{\iter +1}}

\newcommand{\uu}{\bs{u}}
\newcommand{\uk}{\uu_{\iter}}
\newcommand{\ukp}{\uu_{\iter + 1}}

\newcommand{\dcentrk}{d_{\iter}}
\newcommand{\dcentrkp}{d_{\iter +1}}
\newcommand{\lcentrk}{\lambda_{\iter}}
\newcommand{\lcentrkp}{\lambda_{\iter +1}}

\def \algname/{Tracking-ADMM}

\usepackage{soul}

\begin{document}

\begin{frontmatter}
\title{Tracking-ADMM for Distributed Constraint-Coupled Optimization\tnoteref{t1}}
\tnotetext[t1]{This research was supported by the European Research Council (ERC) under the European Union's Horizon 2020 research and innovation programme (grant agreement No 638992 - OPT4SMART) and by the European Commission, under the project UnCoVerCPS, grant number 643921.}
\author[polimi]{Alessandro~Falsone}\ead{alessandro.falsone@polimi.it}
\author[unibo]{Ivano~Notarnicola}\ead{ivano.notarnicola@unibo.it}
\author[unibo]{Giuseppe~Notarstefano}\ead{giuseppe.notarstefano@unibo.it}
\author[polimi]{Maria~Prandini}\ead{maria.prandini@polimi.it}
\address[polimi]{Dip. di Elettronica, Informazione e Bioingegneria,
                    Politecnico di Milano, Via Ponzio 34/5, 20133 Milano, Italy.}
\address[unibo]{Dip. di Ingegneria dell'Energia Elettrica e dell'Informazione ``G. Marconi'',
                    Alma Mater Studiorum Universit\`{a} di Bologna, Viale del Risorgimento 2, 40136 Bologna, Italy.}

\begin{abstract}
We consider constraint-coupled optimization problems in which agents of a network aim to cooperatively minimize the sum of local objective functions subject to individual constraints and a common linear coupling constraint. We propose a novel optimization algorithm that embeds a dynamic average consensus protocol in the parallel Alternating Direction Method of Multipliers (ADMM) to design a fully distributed scheme for the considered set-up. The dynamic average mechanism allows agents to track the time-varying coupling constraint violation (at the current solution estimates). The tracked version of the constraint violation is then used to update local dual variables in a consensus-based scheme mimicking a parallel ADMM step. Under convexity, we prove that all limit points of the agents' primal solution estimates form an optimal solution of the constraint-coupled (primal) problem. The result is proved by means of a Lyapunov-based analysis simultaneously showing consensus of the dual estimates to a dual optimal solution, convergence of the tracking scheme and asymptotic optimality of primal iterates. A numerical study on optimal charging schedule of plug-in electric vehicles corroborates the theoretical results.
\end{abstract}


\end{frontmatter}

\section{Introduction} \label{sec:intro}
Given the increasing penetration of network-connected devices and cloud-computing services in the industry sector, systems are becoming more and more large-scale, composed by different interacting agents, and complex. Efficiently operating these systems is thus becoming more challenging: the computational burden associated to the computation of an efficient strategy for the entire system may be prohibitive for a single processing unit and privacy issues may arise if the different agents are forced to disclose sensitive information to the entity in charge of solving the problem. This challenge has been taken up by the control community, which has directed its efforts in devising distributed algorithms for the resolution of such optimal control problems.

In  this paper we investigate a set-up that is relevant to several network control applications. We consider agents which can communicate only with their neighbors in a network. Each agent knows its local cost function and constraints, that depend on its individual decision variables. Moreover, the agents' decision  variables are coupled by a linear constraint. Each agent knows only how its local decision variables affect the coupling constraint. The agents collectively aim at minimizing the sum of the local cost functions subject to both local and coupling constraints. The presence of this coupling element makes the problem solution challenging, especially in the considered distributed context in which no central authority (communicating with all the agents) is present.

Even though constraint-coupled problems arise more naturally in practical applications the main effort in distributed optimization, until recently, has been devoted to solve optimization problems in which there is no coupling constraint but the agents are required to agree on a \emph{common decision}.
For this set-up, consensus methods based on (sub)gradient iterations and proximal operators (see~\cite{johansson2008subgradient, nedic2009distributed, nedic2010constrained, zanella2011newton, jakovetic2014fast, shi2015extra, nedic2015distributed, margellos2018distributed}), distributed algorithms based on duality (see~\cite{duchi2012dual, zhu2012distributed}), and distributed approaches based on the Alternating Direction Method Multipliers (ADMM) (see~\cite{mota2013dadmm, ling2014decentralized, shi2014linear, jakovetic2015linear, iutzeler2016explicit, makhdoumi2017convergence}) have been proposed.
We refer the reader to~\cite{bertsekas1989parallel, boyd2011distributed} for an in-depth discussion on ADMM.
Recently, distributed gradient schemes for problems with common decision variables have also been combined with a technique known as dynamic average consensus (firstly proposed in~\cite{zhu2010discrete, kia2018tutorial}) to achieve linear convergence rate with constant step-size at the expense of imposing more restrictive assumptions on the cost functions, see~\cite{dilorenzo2016next, varagnolo2016newton, nedic2017achieving, qu2018harnessing, xu2018convergence, xi2018addopt}.
Applying the above mentioned techniques to the constraint-coupled set-up considered in this paper is possible but it would require to define as common decision vector the collection of all local decision vectors. Agents would then have to store and update the entire solution estimate of the problem rather than the portion associated with their decision variables only, and would have all to know the global coupling constraint. These major drawbacks hamper the applicability of those methods to practical applications, thus calling for novel and efficient strategies that takes advantage of the network structure of the problem at hand.

Recently developed approaches tackling constraint-coupled problems directly are mostly based on Lagrangian duality. In~\cite{necoara2015linear} a problem with sparse coupling constraints is considered and a distributed dual gradient method achieving linear convergence rate is proposed. A similar structured problem is considered in~\cite{alghunaim2018dual} where a tailored algorithm exploiting such structure is proved to converge to the optimal solution of the problem.
In~\cite{chang2014distributed,mateos2017distributed}, primal-dual approaches are proposed but they need a diminishing step-size to achieve convergence.
In~\cite{simonetto2016primal, falsone2017dual} distributed dual subgradient algorithms are proposed, in~\cite{chang2016proximal} the dual problem is tackled by means of consensus-ADMM and proximal operators, while an alternative approach based on successive duality steps has been investigated in~\cite{notarnicola2017constraint}; however, all these algorithms typically exhibit a slow convergence rate for the local decision variables. The tracking mechanism has been also employed in~\cite{kia2017distributed} to solve constraint-coupled problems based on an augmented Lagrangian approach for a continuous time setting, but the considered set-up does not allow for nonsmooth costs and local constraints. Finally, in the very recent paper~\cite{zhang2018consensus} an idea similar to the one proposed in this paper is introduced. However, the approach in~\cite{zhang2018consensus} requires the agents to perform multiple communication rounds and the number of rounds has to be carefully tuned to bound a non-vanishing steady state error between the asymptotic algorithmic solution and the optimal one.

The contributions of this paper are as follows. We propose a novel, fully distributed optimization algorithm to solve \emph{constraint-coupled problems} over networks by means of an ADMM-based approach. Differently from distributed ADMM schemes for problems with common decision variables, we design our \algname/ distributed algorithm by embedding a tracking mechanism into the parallel ADMM designed for constraint-coupled problems.

The resulting algorithm enjoys the following appealing features: (i) no parameter tuning is needed, in fact our algorithm works for all the (constant) choices of a penalty parameter and no other coefficients are necessary; (ii) agents solve optimization problems depending on their local (few) decision variables and asymptotically compute only their portion of an optimal (hence feasible) solution to the given problem; (iii) the local estimate of the coupling constraint violation gives each agent a local assessment on the amount of infeasibility, which can be useful, e.g., in designing distributed (receding horizon) control schemes.

The convergence proof of our \algname/ for constraint-coupled problems relies on control systems analysis tools. By explicitly relying on Lyapunov theory for linear systems, we are able to find a constrained solution to a suitable discrete Lyapunov equation that leads to an aggregate descent condition allowing us to prove: (i) the (exact) convergence of the dual variables to a dual optimal solution and (ii) that any limit point of the primal sequences is an optimal solution of the original constraint-coupled problem. This novel approach allows us to derive a clean and elegant proof of the asymptotic optimality of our algorithm.

The rest of the paper is organized as follows. In Section~\ref{sec:setup_preliminaries} we present the problem set-up and in Section~\ref{sec:ADMM_parallel} we revise the ADMM algorithm. In Section~\ref{sec:algorithm} we introduce our novel \algname/ distributed algorithm and analyze its convergence properties, discussing the main steps of the proof. In Section~\ref{sec:simulations} we apply our algorithm on a realistic application related to the optimal charging schedule for a fleet of electric vehicles. In Section~\ref{sec:conclusions} we draw some conclusions and finally, in~\ref{sec:appendix} we report the proofs of all the results stated in the body of the paper.

\paragraph*{Notation}
The vector in $\bR^{n}$ containing all ones is denoted by $\bone_{n}$. The identity matrix of order $n$ is denoted by $I_n$. The Kronecker product is denoted by $\kron$. For a matrix $S$ we write $S\T$ to denote its transpose, $S \succ 0$ to denote that $S$ is positive definite, $\rho(S)$ to denote the spectral radius of $S$, and $\norm{S}$ for its spectral norm. We write $\blkdiag(S_1,\dots,S_n)$ to refer to the block-diagonal matrix with $S_1,\dots,S_n$ as blocks. For a vector $v$, $\norm{v}$ is the Euclidean norm of $v$, and, for any matrix $S \succ 0$, $\norm{v}_S$ is the weighted norm of $v$, i.e., $\norm{v}_S^2 = v\T S v$.

\section{Constraint-Coupled Optimization} \label{sec:setup_preliminaries}
In this section we introduce the optimization set-up and recall some preliminaries about the Alternating Direction Method of Multipliers (ADMM).

\subsection{Optimization Problem and Assumptions} \label{sec:setup}
We consider a system composed of $\Nag$ agents which are willing to cooperatively solve an optimization program formulated over the entire system. Each agent has to set its local decision variables $x_i \in \bR^{n_i}$ so as to minimize the sum of local objective functions $\map{f_i}{\bR^{n_i}}{\bR}$, while satisfying local constraints $X_i \subset \bR^{n_i}$ as well as a linear constraint that couples the decisions of all the agents. Formally, we address the following mathematical program
\begin{align}
    \tag{$\cP$} \label{eq:problem}
    \begin{split}
        \min_{ x_1,\dots,x_\Nag } \;\quad\; &\sum_{i=1}^\Nag f_i(x_i) \\
        \text{subject to:} \quad &\sum_{i=1}^\Nag A_i x_i = b \\
            &x_i \in X_i \qquad i = \toN,
    \end{split}
\end{align}
where $A_i \in \bR^{p \times n_i}$ and $b \in \bR^p$ specify the coupling constraint.

We impose the following regularity conditions on~\ref{eq:problem}.
\begin{assumption}[Convexity and compactness] \label{ass:convexity}
    For all $i = \toN$, the function $f_i$ is convex and the set $X_i$ is convex and compact. \qed
\end{assumption}

Let $\xx = [x_1\T\,\cdots\,x_\Nag\T]\T$, consider a vector $\lambda \in \bR^p$ of Lagrange multipliers and define
\begin{equation} \label{eq:lagrangian}
    L(\xx,\lambda) = \sum_{i=1}^\Nag f_i(x_i) + \lambda\T \bigg(\sum_{i=1}^\Nag A_i x_i - b \bigg)
\end{equation}
the Lagrangian function obtained by dualizing the coupling constraint $\sum_{i=1}^\Nag A_i x_i = b$. Then, the dual problem of~\ref{eq:problem} is
\begin{equation*}
    \max_{\lambda \in \bR^p} \min_{\xx \in X} L(\xx,\lambda) = \max_{\lambda \in \bR^p} \: \: \sum_{i=1}^\Nag \dual_i (\lambda), \tag{$\cD$} \label{eq:dual_problem}
\end{equation*}
where $X = X_1 \times \cdots \times X_\Nag$, and the $i$-th contribution $\dual_i$ is defined as
\begin{equation} \label{eq:local_dual_function}
    \dual_i(\lambda) = \min_{x_i \in X_i} f_i (x_i) + \lambda \T ( A_i x_i - b_i),
\end{equation}
the vectors $b_1,\dots,b_\Nag$ being such that $\sum_{i=1}^\Nag b_i = b$.

The next assumption ensures that~\ref{eq:problem} and~\ref{eq:dual_problem} are well-posed.
\begin{assumption}[Existence of optimal solutions] \label{ass:saddle}
    Problem~\ref{eq:problem} admits an optimal solution $\xs = [x_1^\star\!\T\,\cdots\,x_\Nag^\star\!\T]\T$ and problem~\ref{eq:dual_problem} admits an optimal solution $\ls$. \qed
\end{assumption}

Next, we revise the popular ADMM algorithm which provides an effective way to solve~\ref{eq:problem} by splitting the computation over $\Nag$ agents coordinated by a central unit.

\subsection{The ADMM Parallel Algorithm} \label{sec:ADMM_parallel}
A version of the ADMM algorithm specifically tailored to problem~\ref{eq:problem} is presented in \cite[pag. 254, eq. (4.75)]{bertsekas1989parallel} and is reported here with our notation for the reader's convenience.
Given initial values $x_{i,0} \in X_i$, $d_{0} = (1/\Nag)\sum_{i = 1}^{\Nag} (A_i x_{i,0} - b_i)$, and $\lambda_{0} \in \bR^p$, at each iteration $\iter\ge 0$, a set of $\Nag$ agents and a central unit perform the following two steps. First, all agents, $i = \toN$, compute (in parallel) a minimizer of the following optimization problem
\begin{subequations} \label{eq:parallel_ADMM}
\begin{align}
    \label{eq:centr_primal_update}
    \begin{split}
        &\xikp \in \argmin_{x_i \in X_i} \Big\{ f_i(x_i) + \lambda_\iter \T A_i x_i \\
        &\hspace{9em} +\frac{c}{2} \norm {A_i x_i - A_i \xik + d_\iter }^2 \Big\},
    \end{split}
\end{align}
where $c > 0$ is a (constant) penalty parameter. Then, each agent sends the quantity $A_i \xikp - b_i$ to the central unit, which computes and broadcasts back to all agents the following two quantities
\begin{align}
    &\dcentrkp = \frac{1}{\Nag} \sum_{i=1}^\Nag \left( A_i \xikp - b_i \right) \label{eq:centr_tracker_update} \\
    &\lcentrkp = \lcentrk + c \, \dcentrkp, \label{eq:centr_dual_update}
\end{align}
where the parameter $c$ is the same as in~\eqref{eq:centr_primal_update}. We shall point out that $\dcentrkp$ has no dynamics as it is the \emph{average} of the local contributions $A_i \xikp - b_i$ to the
coupling constraint and it measures the infeasibility of the current tentative solutions $\xikp$, $i = \toN$. Its ``average'' structure is crucial and will be exploited in the design of our \algname/ distributed algorithm.
\end{subequations}

The evolution of~\eqref{eq:parallel_ADMM} is analyzed in \cite[pagg. 254-256]{bertsekas1989parallel} and its convergence property is reported below.
\newtheorem*{proposition*}{Proposition}
\begin{proposition*}[{\cite[Proposition~4.2]{bertsekas1989parallel}}]
    Let Assumptions~\ref{ass:convexity} and~\ref{ass:saddle} hold. Then, any limit point of the primal sequence $\seq{[x_{1,\iter}\T\,\cdots\,x_{\Nag,\iter}\T]\T}$ generated by~\eqref{eq:centr_primal_update}, is an optimal solution of~\ref{eq:problem}, and the dual sequence $\seq{\lcentrk}$, generated by~\eqref{eq:centr_dual_update}, converges to an optimal solution of~\ref{eq:dual_problem}. \qed
\end{proposition*}

Notice that, since all limit points of the primal sequence are optimal, they are necessarily feasible for the coupling constraint. It thus follows that the sequence $\seq{\dcentrk}$ generated by \eqref{eq:centr_tracker_update} converges to zero. Moreover, we shall stress that no requirement on the penalty parameter $c$ is necessary for the convergence result to hold.

Finally, note that the algorithm described by~\eqref{eq:parallel_ADMM} requires a central unit to compute~\eqref{eq:centr_tracker_update} and~\eqref{eq:centr_dual_update}. This hampers the applicability of ADMM to a distributed computation framework, where agents communicate only with neighbors according to a (typically sparse) communication graph. Schemes requiring a central unit are usually referred to as \emph{parallel} schemes. For this reason in the rest of the paper we will refer to~\eqref{eq:parallel_ADMM} as the \emph{parallel ADMM}.


\section{\algname/ Distributed Algorithm} \label{sec:algorithm}
In this section we propose our novel distributed optimization algorithm for constraint-coupled problems. We first introduce the distributed framework and then describe the algorithm along with its convergence properties.

\subsection{Distributed Computation Framework}
Assume that, at each iteration $\iter$, the $\Nag$ agents communicate according to a graph $\cG = (\cV,\cE)$, where $\cV = \until{\Nag}$ is the set of nodes, representing the agents, and $\cE \subseteq \cV \times \cV$ is the set of edges, representing the communication links. The presence of edge $(i,j)$ in $\cE$ models the fact that agent $i$ receives information from agent $j$. We assume that the communication graph does not change across iterations and, consequently, $\cE$ does not depend on the iteration index $\iter$. We denote by $\nbrs_i = \{j \in \cV \mid (i,j) \in \cE \}$ the set of \emph{neighbors} of agent $i$ in $\cG$, assuming that $(i,i) \in \cE$ for all $i = \toN$. We impose the following connectivity property on $\cG$.
\begin{assumption}[Connectivity] \label{ass:connectivity}
    The graph $\cG$ is undirected and connected, i.e., $(i,j)\in\cE$ if and only if $(j,i) \in \cE$ and for every pair of vertices in $\cV$ there exists a path of edges in $\cE$ that connects them. \qed
\end{assumption}

Each edge $(i,j) \in \cE$ has an associated weight $w_{ij}$, which measures how much agent $i$ values the information received by agent $j$. For those $(i,j) \notin \cE$ we set $w_{ij} = 0$, which models the fact that agent $i$ does not receive any information from agent $j$. We impose the following assumption on the network weights.
\begin{assumption}[Balanced information exchange] \label{ass:balance_info}
    For all $i,j = \toN$, $w_{ij} \in [0,1)$ and $w_{ij} = w_{ji}$. Furthermore
    \begin{itemize}
        \item $\sum_{i=1}^\Nag w_{ij} = 1$ for all $j = \toN$,
        \item $\sum_{j=1}^\Nag w_{ij} = 1$ for all $i = \toN$,
    \end{itemize}
    and $w_{ij} > 0$ if and only if $(i,j) \in \cE$. \qed
\end{assumption}
Let $\cW \in \bR^{\Nag \times \Nag}$ be the matrix whose $(i,j)$-th entry is $w_{ij}$, often referred to as the \emph{consensus matrix}. Assumption~\ref{ass:balance_info} translates into requiring $\cW$ to be symmetric and doubly stochastic, i.e., $\cW = \cW\T$ and $\cW \bone_\Nag = \cW\T \bone_\Nag = \bone_\Nag$. We should point out that Assumptions~\ref{ass:connectivity} and~\ref{ass:balance_info} are common in the consensus-based distributed optimization literature, see, e.g., \cite{nedic2009distributed,nedic2010constrained}.

Finally, we impose the following additional assumption on the consensus matrix.
\begin{assumption} \label{ass:W_eigval}
    $\cW$ is positive semidefinite. \qed
\end{assumption}
Note that this assumption is not too restrictive. Indeed, in Section~\ref{sec:W_eigval} we will show that if the agents perform two consecutive communication in the same iteration with any consensus matrix satisfying Assumption~\ref{ass:balance_info}, then this is equivalent to a single communication with a consensus matrix satisfying both Assumptions~\ref{ass:balance_info} and~\ref{ass:W_eigval}.

\subsection{Algorithm Description} \label{sec:algo_description}
In this section, we start from the parallel ADMM and gradually introduce the reader to our proposed algorithm to jointly gain insights about the underlying mechanism and motivate the role of the \emph{consensus and tracking schemes}.

The update~\eqref{eq:centr_dual_update} for~$\lcentrk$ in the parallel ADMM resembles a gradient step aiming at maximizing the cost function of~\ref{eq:dual_problem}. Distributed implementations of gradient-based approaches to solve optimization problems with common decision variables in the form of~\ref{eq:dual_problem} are well known, see e.g.~\cite{nedic2009distributed}.
Typically, each agent $i$ maintain a vector $\lik\in\bR^p$, representing a local version (or copy) of $\lcentrk$, which is iteratively updated according to a consensus-based scheme to force agreement of the local copies.
If, at each iteration, the quantity $\dcentrkp$ were available to all agents, then we could propose the following local update step of $\lik$ for agent $i$
\begin{align} \label{eq:dual_update_description}
    \likp = \sum_{j\in\nbrs_i} w_{ij} \, \ljk + c \, \dcentrkp,
\end{align}
for all $i = \toN$.

However, update~\eqref{eq:dual_update_description} cannot be implemented in a fully distributed scheme since $\dcentrkp$ is not locally available and should be computed by a central unit (cf.~\eqref{eq:centr_tracker_update}). In the same spirit of $\lik$, a distributed counterpart for \eqref{eq:centr_tracker_update} can be obtained by equipping each agent $i$ with a local auxiliary quantity $\dik \in \bR^p$, which serves as a local estimate of $\dcentrk$. Since $\dcentrk$ is the average of $A_i \xik - b_i$, $i = \toN$ (with $\sum_{i=1}^\Nag b_i = b$), we propose to update $\dik$ according to a (distributed) dynamic average consensus mechanism, see \cite{zhu2010discrete,kia2018tutorial}:
\begin{equation} \label{eq:tracker_update_description}
    \dikp = \sum_{j\in\nbrs_i} w_{ij} \, \djk + (A_i \xikp - b_i) - (A_i \xik - b_i),
\end{equation}
initialized with $d_{i,0} = A_i x_{i,0} - b_i$, for all $i = \toN$.

In this way, $\dik$ acts as a \emph{distributed tracker} of the (time-varying) signal $(1/\Nag) \sum_{i=1}^\Nag \left(A_i \xik - b_i\right)$. Using $\dikp$ in place of $\dcentrkp$ in~\eqref{eq:dual_update_description} makes the update of $\lik$ fully distributed.

Clearly, since we propose to replace the centralized quantities $\lcentrk$ and $\dcentrk$ with their local counterparts, the local minimization in~\eqref{eq:centr_primal_update} has to be adjusted accordingly.
The \algname/ is formally summarized in Algorithm~\ref{alg:DT-ADMM} from the perspective of agent $i$.
Specifically, the adapted local minimization to compute $\xikp$ is shown in~Step~\ref{step:primal_update}, where $\lcentrk$ and $\dcentrk$ in the original centralized update~\eqref{eq:centr_primal_update}, are replaced by the local averages $\ellik$ and $\deltaik$, respectively (cf. Steps~\ref{step:tracker_consensus} and~\ref{step:dual_consensus}).

\begin{algorithm}[t]
    \begin{spacing}{\algspc}
    \begin{algorithmic}[1]
        \STATE \bf{Initialization}
        \STATE ~~~~$x_{i,0} \in X_i$ \label{step:primal_init}
        \STATE ~~~~$d_{i,0} = A_i x_{i,0} - b_i$ \label{step:tracker_init}
        \STATE ~~~~$\lambda_{i,0} \in \bR^p$ \label{step:dual_init}
        \STATE \bf{Repeat until convergence}
        \STATE ~~~~$\deltaik = \sum_{j\in \nbrs_i} w_{ij} \, \djk$ \label{step:tracker_consensus}
        \STATE ~~~~$\ellik = \sum_{j\in \nbrs_i} w_{ij} \, \ljk$ \label{step:dual_consensus}
        \STATE ~~~~$\xikp \in \argmin\limits_{x_i \in X_i} \Big\{ f_i(x_i) + \ellik\T A_i x_i$ \\ \vspace{-0.3cm}
            \begin{flushright}$+\dfrac{c}{2} \norm {A_i x_i - A_i \xik + \deltaik }^2 \Big\}$\end{flushright}
            \label{step:primal_update}\vspace{0.0cm}
        \STATE ~~~~$\dikp = \deltaik + A_i \xikp - A_i \xik $ \label{step:tracker_update}
        \STATE ~~~~$\likp = \ellik + c \, \dikp$ \label{step:dual_update}
        \STATE ~~~~$\iter \leftarrow \iter+1$
    \end{algorithmic}
    \end{spacing}
    \caption{\algname/}
    \label{alg:DT-ADMM}
\end{algorithm}

Some remarks are in order. First, we shall stress that all update steps (cf. Steps~\ref{step:primal_update}-\ref{step:dual_update}) are fully distributed, as they only use quantities that are either locally known to agent $i$ or collected by agent $i$ via neighboring communications (cf. Steps~\ref{step:tracker_consensus} and~\ref{step:dual_consensus}). Moreover, the minimization in Step~\ref{step:primal_update} is always well defined in view of Assumption~\ref{ass:convexity}.

While the initialization of $\xik$ and $\lik$ can be arbitrary, the correct initialization of $\dik$ as per Step~\ref{step:tracker_init} is crucial, consistently with other tracking-based approaches as, e.g., the ones mentioned in the introduction (cf.~\cite{dilorenzo2016next, varagnolo2016newton, nedic2017achieving, qu2018harnessing, xu2018convergence, xi2018addopt}).
Sensible values for initializing Algorithm~\ref{alg:DT-ADMM} are given by $x_{i,0} \in \argmin_{x_i \in X_i} f_i(x_i)$, $d_{i,0} = A_i x_{i,0} - b/\Nag$, and $\lambda_{i,0} = 0$.

Finally, the parameter $c > 0$ in Step~\ref{step:primal_update} and~\ref{step:dual_update} is constant and is similar to the step-size of a gradient-based method, but, differently from gradient-like approaches, its value can be arbitrary.

\subsection{Algorithm Analysis}
In this subsection we analyze the proposed \algname/ algorithm and state its convergence properties. To ease the exposition, we state and discuss the main steps of the convergence analysis, deferring their proofs to \ref{sec:appendix}.

\subsubsection{Aggregate Reformulation of \algname/} \label{sec:global_reform}
First of all we shall rewrite the execution of Algorithm~\ref{alg:DT-ADMM} by all agents in a compact form.

Recalling that $\xx = [x_1\T\,\cdots\,x_\Nag\T]\T$, let us consistently denote with bold symbols the vectors collecting the corresponding quantity of all agents, i.e.,
\begin{align*}
    \xk &= [x_{1,\iter}\T\,\cdots\,x_{\Nag,\iter}\T]\T, \\
    \dk &= [d_{1,\iter}\T\,\cdots\,d_{\Nag,\iter}\T]\T,
        &\deltak &= [\delta_{1,\iter}\T\,\cdots\,\delta_{\Nag,\iter}\T]\T, \\
    \lk &= [\lambda_{1,\iter}\T\,\cdots\,\lambda_{\Nag,\iter}\T]\T,
        &\ellk &= [\ell_{1,\iter}\T\,\cdots\,\ell_{\Nag,\iter}\T]\T.
\end{align*}

Then, the evolution of Algorithm~\ref{alg:DT-ADMM} over the whole multi-agent network can be compactly written as
\begin{subequations} \label{eq:DT-ADMM}
    \begin{align}
        \begin{split} \label{eq:all_primal_update}
        &\xkp \in \argmin_{\xx \in X} \Big\{ f(\xx) + (W \lk)\T A_d \xx \\
            &\hspace{8em} +\frac{c}{2} \norm {A_d \xx - A_d \xk + W \dk }^2 \Big\}
        \end{split}
        \\
        &\dkp = W \dk + A_d \xkp - A_d \xk \label{eq:all_tracker_update} \\
        &\lkp = W \lk + c \, \dkp, \label{eq:all_dual_update}
    \end{align}
\end{subequations}
where $f(\xx) = \sum_{i=1}^\Nag f_i(x_i)$, $A_d = \blkdiag(A_1,\dots,A_\Nag)$, $X = X_1 \times \cdots \times X_\Nag$, and $W = \cW \kron I_p$. Note that, in~\eqref{eq:DT-ADMM}, we used $\deltak = W \dk$ and $\ellk = W \lk$, which represent the network-wide formulations of Steps~\ref{step:tracker_consensus} and~\ref{step:dual_consensus} of Algorithm~\ref{alg:DT-ADMM}, respectively.

As can be seen from~\eqref{eq:DT-ADMM}, \algname/ describes a dynamical system composed of two parts: a nonlinear dynamics given by the optimization step (cf.~\eqref{eq:all_primal_update}) and a linear dynamics given by a twofold consensus-based step (cf.~\eqref{eq:all_tracker_update} and~\eqref{eq:all_dual_update}). In the forthcoming analysis, we will investigate the structural properties of these two parts separately. Then we show how to combine them to prove convergence of the proposed algorithm to an optimal solution of the original problem~\ref{eq:problem}.

\subsubsection{Discussion on Assumption~\ref{ass:W_eigval}} \label{sec:W_eigval}
Before starting with the convergence analysis, we now show that Assumption~\ref{ass:W_eigval} is not restrictive.

Suppose to substitute Steps~\ref{step:tracker_consensus} and~\ref{step:dual_consensus} with
\begin{equation*}
    \deltaik' = \textstyle \sum_{j\in \nbrs_i} w_{ij} \, \djk
    \quad\text{and}\quad
    \ellik' = \textstyle \sum_{j\in \nbrs_i} w_{ij} \, \ljk
\end{equation*}
and then perform the updates
\begin{equation*}
    \deltaik = \textstyle \sum_{j\in \nbrs_i} w_{ij} \, \deltajk'
    \quad\text{and}\quad
    \ellik = \textstyle \sum_{j\in \nbrs_i} w_{ij} \, \elljk'
\end{equation*}
before Step~\ref{step:primal_update}. These two steps can then be compactly written as
\begin{align*}
    \deltak &= W^2 \dk = (\cW \kron I_p)^2 \dk = (\cW^2 \kron I_p) \dk, \\
    \ellk &= W^2 \lk = (\cW \kron I_p)^2 \lk = (\cW^2 \kron I_p) \lk,
\end{align*}
which are equivalent to a single update round with $\cW^2$ in place of $\cW$. Since, under Assumption~\ref{ass:balance_info}, $\cW$ is symmetric and doubly stochastic, then $\cW^2$ is symmetric, doubly stochastic, and positive semidefinite, thus satisfying Assumption~\ref{ass:W_eigval}. Letting $\cE^2$ be the set of edges $(i,j)$ such that there exists a $\kappa \in \cV$ for which $(i,\kappa),(\kappa,j) \in \cE$, it is easy to show that the graph $\cG^2 = (\cV,\cE^2)$ satisfies Assumption~\ref{ass:connectivity} and $\cW^2$ satisfies Assumption~\ref{ass:balance_info} with $\cE^2$ in place of $\cE$.

Therefore, Assumption~\ref{ass:W_eigval} can be easily satisfied by simply performing two update rounds rather than one. We choose not to alter Algorithm~\ref{alg:DT-ADMM} to keep the notation light and avoid writing $\cW^2$ in place of $\cW$ everywhere.

\subsubsection{Network Average Quantities}
We start the analysis of \algname/ by stating some important properties regarding the network averages of the local variables $\dik$ and $\lik$, which are defined as
\begin{align}
    \dkavg = \frac{1}{\Nag} \sum_{i = 1}^\Nag \dik
    \quad\text{and}\quad
    \lkavg = \frac{1}{\Nag} \sum_{i = 1}^\Nag \lik,
    \label{eq:avg_definition}
\end{align}
respectively. The following lemmas highlight key features of $\dkavg$ and $\lkavg$.

\begin{lemma}[Tracking property] \label{lem:tracking_property}
    Under Assumption~\ref{ass:balance_info} it holds that
    \begin{align}
        \dkavg &= \frac{1}{\Nag} \left( \sum_{i=1}^\Nag A_i \xik - b \right)\!,
        \label{eq:avg_tracker_prop}
    \end{align}
    for all $\iter \geq 0$. \qed
\end{lemma}

\begin{lemma}[Average dual update] \label{lem:avg_dual_prop}
    Under Assumption~\ref{ass:balance_info}, it holds that
    \begin{align}
        \lkpavg & = \lkavg + c \, \dkpavg,
        \label{eq:avg_dual_prop}
    \end{align}
    for all $\iter \geq 0$. \qed
\end{lemma}

As can be seen by a formal comparison, \eqref{eq:avg_tracker_prop} and~\eqref{eq:avg_dual_prop} mimic the update step~\eqref{eq:centr_tracker_update} and~\eqref{eq:centr_dual_update} of the parallel ADMM. Therefore, by virtue of Lemmas~\ref{lem:tracking_property} and~\ref{lem:avg_dual_prop}, running Algorithm~\ref{alg:DT-ADMM} results in the agents network behaving, on average, like the parallel ADMM.

\subsubsection{Consensus Errors as Linear Dynamical Systems}
The introduction of $\dkavg$ and $\lkavg$ is also useful to analyze the behavior of $\dik$ and $\lik$. Specifically, as usually done in consensus-based optimization algorithms, we shall study the consensus error, i.e., distance between the agents' local quantities $\dik$ and $\lik$ and their respective network average $\dkavg$ and $\lkavg$.

To ease the notation, let us introduce the following vectors containing $\Nag$ copies of the respective average quantities
\begin{align*}
    &\ddkavg = \bone_{\Nag} \kron \dkavg, \\
    &\llkavg = \bone_{\Nag} \kron \lkavg,
\end{align*}
which are useful when comparing the aggregate vectors $\dk$ and $\lk$ with their averages $\dkavg$ and $\lkavg$. Then, define the consensus errors as
\begin{align*}
    \edk & = \dk - \ddkavg, \\
    \elk & = \lk - \llkavg,
\end{align*}
and the auxiliary sequence
\begin{equation} \label{eq:z_def}
    \zk = A_d \xk - \ddkavg.
\end{equation}

From the expression of $\dkp$ in~\eqref{eq:all_tracker_update}, we can write
\begin{align}
    \edkp &= \dkp - \ddkpavg \nonumber \\
        & = W \dk + A_d \xkp - A_d \xk - \ddkpavg \pm \ddkavg \nonumber \\
        & \eqlabel{(a)} W \dk - \ddkavg + [\zkp - \zk] \nonumber \\
        & \eqlabel{(b)} W \edk + [\zkp - \zk], \label{eq:stable_ed_dyn}
\end{align}
where in (a) we used the definition of $\zkp$ and $\zk$ (cf. \eqref{eq:z_def}) and in (b) we exploited the following
\begin{align*}
    W \ddkavg &= (\cW \kron I_p)(\bone_\Nag \kron \dkavg) \\
        & \eqlabel{(c)} (\cW \bone_\Nag \kron I_p \dkavg) \\
        & \eqlabel{(d)} (\bone_\Nag \kron \dkavg) = \ddkavg,
\end{align*}
which holds due to the mixed-product property of the Kronecker product in (c) and the doubly stochasticity of $\cW$ under Assumption~\ref{ass:balance_info} in (d).

Finally, letting $\ones$ be the $\Nag\times \Nag$ matrix with all entries equal to $1/\Nag$ and noticing that
\begin{align*}
    (\ones \kron I_p) \edk &= \left( \frac{1}{\Nag}\bone_\Nag \bone_\Nag\T \kron I_p \right) [\dk - \ddkavg] \\
        &= \bone_\Nag \kron \frac{1}{\Nag}\sum_{i=1}^\Nag [\dik - \dkavg] \\
        &= \bone_\Nag \kron \frac{1}{\Nag} [\Nag \dkavg - \Nag \dkavg] = 0,
\end{align*}
we can subtract the quantity $(\ones \kron I_p) \edk = 0$ from the right hand side of~\eqref{eq:stable_ed_dyn} to get
\begin{equation} \label{eq:as_stable_ed_dyn}
    \edkp = \tW \edk + [\zkp - \zk],
\end{equation}
where we set $\tW = W - (\ones \kron I_p)$.

With an analogous reasoning, we can also show that
\begin{align}
    \elkp &= \lkp - \llkpavg \nonumber \\
        & \eqlabel{(a)} W \lk + c\,\dkp - (\llkavg + c\,\ddkpavg) \nonumber \\
        &= W \elk + c\,\edkp \nonumber \\
        &= \tW \elk + c\,\edkp, \label{eq:as_stable_el_dyn}
\end{align}
where, in (a), we used~\eqref{eq:all_dual_update} and~\eqref{eq:avg_dual_prop} pre-multiplied by $\bone_\Nag \kron$.

Under Assumptions~\ref{ass:connectivity} and~\ref{ass:balance_info}, we have that $\norm{\cW - \ones} < 1$, then
\begin{equation*}
    \rho(\tW) \leq \norm{\tW} = \norm{(\cW - \ones) \kron I_p} = \norm{\cW - \ones} \norm{I_p} < 1.
\end{equation*}
Thus, all eigenvalues of $\tW$ lie within the unit circle. This implies that the dynamical systems in~\eqref{eq:as_stable_ed_dyn} and~\eqref{eq:as_stable_el_dyn}, describing the evolution of the consensus errors, are asymptotically stable. The importance of this fact is twofold: i) \eqref{eq:as_stable_ed_dyn} and~\eqref{eq:as_stable_el_dyn} enjoys the input-to-state stability property: if the sequence $\seq{\zkp - \zk}$ is bounded then also $\seq{\edk}$ and $\seq{\elk}$ are bounded; ii) if $\zkp - \zk$ vanishes, then also $\edk$ and $\elk$ do.

Point i) is formalized in the following lemma.
\begin{lemma}[Bounded sequences] \label{lem:bounded_errors}
    Under Assumptions~\ref{ass:convexity},~\ref{ass:connectivity} and~\ref{ass:balance_info} we have that
    \begin{enumerate}[(i)]
        \item the sequences $\seq{\xk}$ and $\seq{\zk}$ are bounded; \label{item:x_z_bounded}
        \item the sequences $\seq{\edk}$ and $\seq{\elk}$ are bounded. \qed \label{item:ed_el_bounded}
    \end{enumerate}
\end{lemma}

The results in Lemmas~\ref{lem:tracking_property}, \ref{lem:avg_dual_prop}, and~\ref{lem:bounded_errors} can be interpreted as follows. The sequences $\seq{\dk}$ and $\seq{\lk}$ generated by \eqref{eq:all_tracker_update} and~\eqref{eq:all_dual_update} evolve in a (possibly large) neighborhood of their averages $\dkavg$ and $\lkavg$, which, in turn, evolve according to a parallel ADMM iteration.

\subsubsection{Optimality Characterization}
We now focus our attention to the analysis of the (nonlinear) update step in \eqref{eq:all_primal_update}. The following result represents the distributed counterpart of the key inequality to prove convergence of the parallel ADMM.
\begin{proposition}[Local optimality] \label{prop:optimality_descent}
    Under Assumptions~\ref{ass:convexity} and~\ref{ass:saddle} we have that
    \begin{equation} \label{eq:optimality_descent}
        \begin{multlined}
            \norm{\llkpavg - \lls}^2 + 2 c [\zkp - A_d \xs]\T \elkp \\
                \leq \norm{\llkavg - \lls}^2 - \norm{\llkpavg - \llkavg}^2,
        \end{multlined}
    \end{equation}
    for any optimal solution pair $(\xs,\ls)$ for~\ref{eq:problem} and~\ref{eq:dual_problem}, where we set $\lls = \bone_{\Nag} \kron \ls$. \qed
\end{proposition}

If we ignore the cross term $2 c [\zkp - A_d \xs]\T \elkp$, then~\eqref{eq:optimality_descent} tells us that as long as $\llkavg$ changes from one iteration to the next one, the network of agents is, on average, getting closer to the optimal solution of~\ref{eq:dual_problem}. The additional cross term depends on the consensus error $\elkp = \lkp - \llkpavg$ and is clearly a distinctive feature of \algname/ (as opposed to the parallel ADMM) arising from its distributed nature. In the parallel ADMM, a condition similar to \eqref{eq:optimality_descent} can be derived and leveraged to prove the algorithm convergence. However, as expected, this is not the case here:~\eqref{eq:optimality_descent} is not sufficient to prove convergence of the proposed scheme and we are required to study the interplay between the minimization step in~\eqref{eq:all_primal_update} and the two consensus steps in~\eqref{eq:all_tracker_update} and~\eqref{eq:all_dual_update}.

\subsubsection{\algname/ Convergence}
The convergence of Algorithm~\ref{alg:DT-ADMM} is divided into two results. The first theorem proves that the two consensus mechanism are successful, i.e., agents eventually agree on a common value for the tracker and for the dual variables (cf. points~\ref{item:thm_tracker_consensus} and~\ref{item:thm_dual_consensus}), and the algorithm converges to a solution that is feasible for~\ref{eq:problem} (cf. point~\ref{item:thm_feasibility}); the second theorem proves that the primal and dual solution obtained by Algorithm~\ref{alg:DT-ADMM} are optimal for~\ref{eq:problem} and~\ref{eq:dual_problem}, respectively.

\begin{theorem}[Convergence] \label{thm:convergence}
    Under Assumptions~\ref{ass:convexity}-\ref{ass:W_eigval}, the sequences generated by \algname/ satisfy:
    \begin{enumerate}[(i)]
        \item $\lim_{\iter \to \infty} \norm{\edk} = 0$, \label{item:thm_tracker_consensus}
        \item $\lim_{\iter \to \infty} \norm{\elk} = 0$, \label{item:thm_dual_consensus}
        \item $\lim_{\iter \to \infty} \norm{\dkavg} = 0$, \label{item:thm_feasibility}
        \item $\seq{ \norm{\llkavg - \lls}^2 + c^2 \norm{\zk - A_d \xs}^2 }$ is convergent, \label{item:thm_convergence}
    \end{enumerate}
    for any optimal solution pair $(\xs,\ls)$ for~\ref{eq:problem} and~\ref{eq:dual_problem}, with $\lls = \bone_{\Nag} \kron \ls$. \qed
\end{theorem}

The main idea behind the proof of Theorem~\ref{thm:convergence} is the following. We start by considering~\eqref{eq:as_stable_ed_dyn} and~\eqref{eq:as_stable_el_dyn} together as a single dynamical system for the consensus updates (the consensus system).
Since the consensus system is asymptotically stable, then we can build a positive-definite quadratic Lyapunov function which monotonically decreases along its free evolution.
However, the consensus system is \emph{not} autonomous, thus the inequality describing the variation of such Lyapunov function across iterations contains terms that are not defined in sign, so that it is no longer necessarily decreasing. The idea is to properly combine such inequality with~\eqref{eq:optimality_descent} of Proposition~\ref{prop:optimality_descent} so as to ``balance out'' the terms that are not defined in sign. In this way, we obtain an aggregate descent condition that allows us to prove the asymptotic stability of the overall (nonlinear) dynamical system modeling our \algname/ algorithm.

The reader should note that point~\ref{item:thm_tracker_consensus} of Theorem~\ref{thm:convergence} implies that, for all $i = \toN$, $\dik$ tends to $\dkavg$ as $\iter \to \infty$. In light of~\eqref{eq:avg_tracker_prop}, this ensures that each agent is able to locally assess the amount of infeasibility of the current primal iterates for the coupling constraint. This feature is of great importance in those applications where feasibility up to a given tolerance is sufficient.

Finally, the following result shows that by running Algorithm~\ref{alg:DT-ADMM} agents are able to compute optimal solutions to \ref{eq:problem} and~\ref{eq:dual_problem} in a fully distributed way.
\begin{theorem}[Optimality] \label{thm:main_result}
    Under Assumptions~\ref{ass:convexity}-\ref{ass:W_eigval}, the sequences generated by \algname/ are such that
    \begin{enumerate}[(i)]
        \item any limit point of the primal sequence $\seq{\xk}$ is an optimal solution $\xs$ of~\ref{eq:problem};\label{eq:thm_primal_opt}
        \item each dual sequence $\seq{\lik}$, $i = \toN$, converges to the same optimal solution $\ls$ of~\ref{eq:dual_problem}. \qed \label{eq:thm_dual_opt}
    \end{enumerate}
\end{theorem}

We point out that Theorem~\ref{thm:main_result} guarantees that the sequences $\seq{\xik}$ are asymptotically optimal, hence feasible for both local and the coupling constraints. Notice also that Algorithm~\ref{alg:DT-ADMM} works with (not necessarily strictly) convex primal cost functions and does not require any primal recovery procedure.

\section{Numerical Study} \label{sec:simulations}
To demonstrate the effectiveness of our approach, we test \algname/ on a modified version of the Plug-in Electric Vehicles Optimal Charging Schedule problem described in \cite{vujanic2016decomposition}. A fleet of $\Nag$ electric vehicles has to select an overnight charging schedule so as to minimize the sum of the electricity costs for charging the internal battery of each vehicle while complying with vehicle-level constraints (user requirements on the final state of charge and battery physical limitations) and a network-wide constraint (maximum power that the grid can deliver). For simplicity we focus on the ``only charging'' case of the charging problem in \cite{vujanic2016decomposition}.
Differently from \cite{vujanic2016decomposition}, where the vehicles can decide only whether to charge or not their internal battery at a fixed charging rate, here we allow each vehicle to optimize also the charging rate.
The resulting optimization program can then be formalized as follows
\begin{align} \label{eq:PEVs}
\begin{split}
    \min_{ \xi_1,\dots,\xi_\Nag } \;\quad\; &\sum_{i=1}^\Nag \gamma_i\T \xi_i \\
    \text{subject to:} \quad &\sum_{i=1}^\Nag \tilde{A}_i \xi_i \leq b \\
        &\xi_i \in \Xi_i \qquad i = \toN,
\end{split}
\end{align}
where $X_i$ are bounded polyhedral sets. We refer the reader to \cite{vujanic2016decomposition} for the precise formulations of the quantities in~\eqref{eq:PEVs}.

The coupling constraints in \eqref{eq:PEVs} is an inequality, which makes \eqref{eq:PEVs} not fitting the formulation in~\ref{eq:problem}. However, we can easily turn the coupling constraint into an equality by introducing $p$ additional slack variables for each vehicle. Problem~\eqref{eq:PEVs} thus becomes
\begin{align} \label{eq:PEVs_slack}
\begin{split}
    \min_{ \xi_1,\dots,\xi_\Nag } \;\quad\; &\sum_{i=1}^\Nag \gamma_i\T \xi_i \\
    \text{subject to:} \quad &\sum_{i=1}^\Nag ( \tilde{A}_i \xi_i + s_i ) = b \\
        \begin{split}
            &\xi_i \in \Xi_i \\
            &0\leq s_i \leq \bar{s}
        \end{split}
        \qquad i = \toN,
\end{split}
\end{align}
which now fits~\ref{eq:problem} with the identifications $x_i = [\xi_i\T \; s_i\T]\T$, $f_i(x_i) = \gamma_i\T \xi_i$, $A_i = [\tilde{A}_i \; I_p]$, and $X_i = \Xi_i \times S_i$ where $S_i = \{s_i: 0\leq s_i \leq \bar{s} \}$, for some $\bar{s} \ge b$, for all $i = \toN$.

In our simulation we considered $\Nag = 100$ vehicles. Each vehicle $i$ has 24 decision variables $\xi_i$ plus 24 additional slack variables $s_i$ ($n_i = 48$) and the local constraint set $X_i$ is defined by 245 inequalities. The number of coupling constraint, and consequently the number of Lagrange multipliers, is $p = 24$. The communication network and the weight matrix $\cW$ have been generated at random so as to satisfy Assumptions~\ref{ass:connectivity}, \ref{ass:balance_info}, and~\ref{ass:W_eigval}.

We run our \algname/ for 200 iterations, initialized as suggested in Section~\ref{sec:algo_description}, setting $c = 10^{-4}$.

\begin{figure}[t]
    \centering
    \includegraphics[width=\columnwidth]{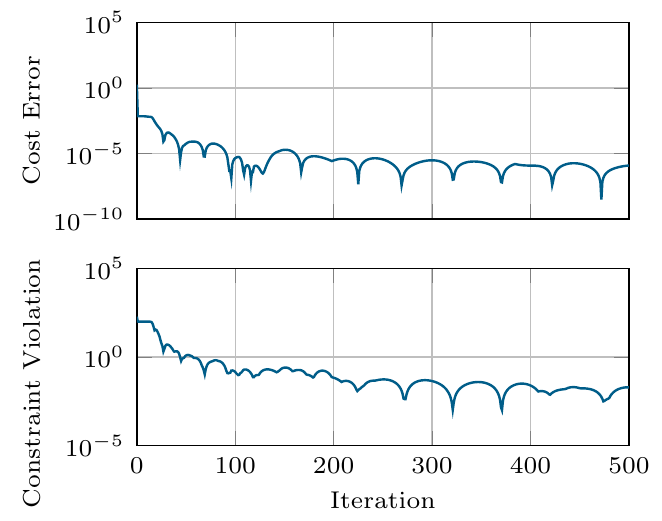}
    \caption{Difference between the cost achieved by the primal tentative solution and the optimal cost (upper plot) and maximum violation of the joint constraints (lower plot), across iterations.}
    \label{fig:cost_constr}
\end{figure}

In Figure~\ref{fig:cost_constr} we report the behavior across iterations of the difference $\sum_{i=1}^{\Nag} f_i(\xik) - f^\star$ between the value of the cost function achieved by the primal tentative solution $\xik$ and the optimal cost $f^\star$ computed by a centralized solver (upper plot), and the infinite norm $\norm{\sum_{i=1}^{\Nag} A_i \xik - b}_\infty$ of the joint constraints violation (lower plot).
As it can be seen from the picture, the proposed algorithm reaches feasibility and optimality. Note that, differently from the approaches based on dual decomposition, the primal iterates generated by Algorithm~\ref{alg:DT-ADMM} achieve feasibility of the coupling constraint without requiring any recovery procedure (which is known to slow down the algorithm convergence).

\begin{figure}[t]
    \centering
    \includegraphics[width=\columnwidth]{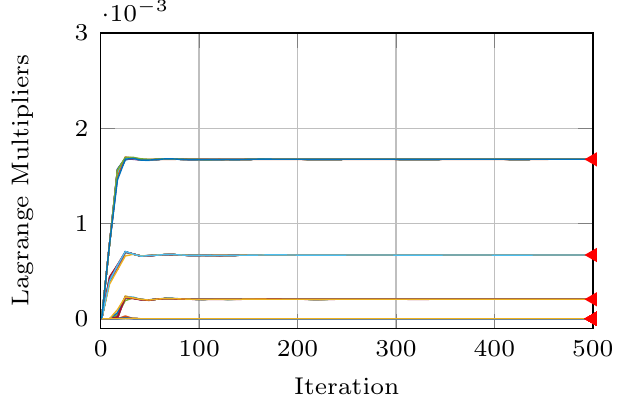}
    \caption{Evolution of Lagrange multipliers across iterations. Red triangles denote the optimal dual solution.}
    \label{fig:lambdas}
\end{figure}

In Figure~\ref{fig:lambdas} we also plot the evolution of the agents local estimates $\seq{\lik}$ of the optimal Lagrange multipliers $\ls$ (red triangles). From the picture the reader can easily see that local estimates reach consensus and optimality.

\section{Conclusions} \label{sec:conclusions}
In this paper we have proposed a novel distributed method to solve constraint-coupled convex optimization problems in which the sum of local cost functions needs to be minimized while satisfying both individual constraints (involving one component of the decision variables) and a common linear coupling constraint (involving all the components). The distributed algorithm combines the (parallel) ADMM algorithm tailored for this class of optimization problems with a dynamic tracking mechanism. We proved that each agent asymptotically computes an optimal dual solution and its portion of an optimal solution to the target (primal) problem. Moreover, the tracking scheme allows agents to obtain a local measure of the coupling-constraint violation. Numerical computations corroborated the theoretical results.

\bibliography{DT-ADMM}

\begin{thebibliography}{10}
\expandafter\ifx\csname url\endcsname\relax
  \def\url#1{\texttt{#1}}\fi
\expandafter\ifx\csname urlprefix\endcsname\relax\def\urlprefix{URL }\fi
\expandafter\ifx\csname href\endcsname\relax
  \def\href#1#2{#2} \def\path#1{#1}\fi

\bibitem{johansson2008subgradient}
B.~Johansson, T.~Keviczky, M.~Johansson, K.~H. Johansson, Subgradient methods
  and consensus algorithms for solving convex optimization problems, in: {IEEE}
  Conference on Decision and Control {(CDC)}, 2008, pp. 4185--4190.

\bibitem{nedic2009distributed}
A.~Nedi{\'c}, A.~Ozdaglar, Distributed subgradient methods for multi-agent
  optimization, {IEEE} Trans. on Automatic Control 54~(1) (2009) 48--61.

\bibitem{nedic2010constrained}
A.~Nedi{\'c}, A.~Ozdaglar, P.~A. Parrilo, Constrained consensus and
  optimization in multi-agent networks, {IEEE} Trans. on Automatic Control
  55~(4).

\bibitem{zanella2011newton}
F.~Zanella, D.~Varagnolo, A.~Cenedese, G.~Pillonetto, L.~Schenato,
  {N}ewton-{R}aphson consensus for distributed convex optimization, in: {IEEE}
  Conference on Decision and Control and European Control Conference
  {(CDC-ECC)}, 2011, pp. 5917--5922.

\bibitem{jakovetic2014fast}
D.~Jakoveti{\'c}, J.~Xavier, J.~M. Moura, Fast distributed gradient methods,
  IEEE Transactions on Automatic Control 59~(5) (2014) 1131--1146.

\bibitem{shi2015extra}
W.~Shi, Q.~Ling, G.~Wu, W.~Yin, Extra: An exact first-order algorithm for
  decentralized consensus optimization, SIAM Journal on Optimization 25~(2)
  (2015) 944--966.

\bibitem{nedic2015distributed}
A.~Nedi{\'c}, A.~Olshevsky, Distributed optimization over time-varying directed
  graphs, IEEE Transactions on Automatic Control 60~(3) (2015) 601--615.

\bibitem{margellos2018distributed}
K.~Margellos, A.~Falsone, S.~Garatti, M.~Prandini, Distributed constrained
  optimization and consensus in uncertain networks via proximal minimization,
  IEEE Transactions on Automatic Control 63~(5) (2018) 1372--1387.

\bibitem{duchi2012dual}
J.~C. Duchi, A.~Agarwal, M.~J. Wainwright, Dual averaging for distributed
  optimization: Convergence analysis and network scaling, IEEE Trans. on Autom.
  Control 57~(3) (2012) 592--606.

\bibitem{zhu2012distributed}
M.~Zhu, S.~Mart{\'i}nez, On distributed convex optimization under inequality
  and equality constraints, IEEE Transactions on Automatic Control 57~(1)
  (2012) 151--164.

\bibitem{mota2013dadmm}
J.~F. Mota, J.~M. Xavier, P.~M. Aguiar, M.~P{\"u}schel, {D}-{ADMM}: A
  communication-efficient distributed algorithm for separable optimization,
  IEEE Trans. on Signal Processing 61~(10) (2013) 2718--2723.

\bibitem{ling2014decentralized}
Q.~Ling, A.~Ribeiro, Decentralized dynamic optimization through the alternating
  direction method of multipliers, {IEEE} Trans. on Signal Processing 5~(62)
  (2014) 1185--1197.

\bibitem{shi2014linear}
W.~Shi, Q.~Ling, K.~Yuan, G.~Wu, W.~Yin, On the linear convergence of the
  {ADMM} in decentralized consensus optimization, {IEEE} Trans. on Signal
  Processing 62~(7) (2014) 1750--1761.

\bibitem{jakovetic2015linear}
D.~Jakoveti{\'c}, J.~M. Moura, J.~M. Xavier, Linear convergence rate of a class
  of distributed augmented {L}agrangian algorithms, IEEE Trans. on Automatic
  Control 60~(4) (2015) 922--936.

\bibitem{iutzeler2016explicit}
F.~Iutzeler, P.~Bianchi, P.~Ciblat, W.~Hachem, Explicit convergence rate of a
  distributed alternating direction method of multipliers, {IEEE} Trans. on
  Autom. Control 61~(4) (2016) 892--904.

\bibitem{makhdoumi2017convergence}
A.~Makhdoumi, A.~Ozdaglar, Convergence rate of distributed {ADMM} over
  networks, {IEEE} Trans. on Automatic Control 62~(10) (2017) 5082--5095.

\bibitem{bertsekas1989parallel}
D.~P. Bertsekas, J.~N. Tsitsiklis, Parallel and distributed computation:
  numerical methods, Vol.~23, Prentice hall Englewood Cliffs, NJ, 1989.

\bibitem{boyd2011distributed}
S.~Boyd, N.~Parikh, E.~Chu, B.~Peleato, J.~Eckstein, Distributed optimization
  and statistical learning via the alternating direction method of multipliers,
  Found. and Trends{\textregistered} in Machine learning 3~(1) (2011) 1--122.

\bibitem{zhu2010discrete}
M.~Zhu, S.~Mart{\'i}nez, Discrete-time dynamic average consensus, Automatica
  46~(2) (2010) 322--329.

\bibitem{kia2018tutorial}
S.~S. Kia, B.~Van~Scoy, J.~Cortes, R.~A. Freeman, K.~M. Lynch, S.~Martinez,
  Tutorial on dynamic average consensus: The problem, its applications, and the
  algorithms, IEEE Control Systems Magazine 39~(3) (2019) 40--72.

\bibitem{dilorenzo2016next}
P.~Di~Lorenzo, G.~Scutari, {NEXT}: In-network nonconvex optimization, {IEEE}
  Trans. on Signal and Information Process. over Networks 2~(2) (2016)
  120--136.

\bibitem{varagnolo2016newton}
D.~Varagnolo, F.~Zanella, A.~Cenedese, G.~Pillonetto, L.~Schenato,
  {N}ewton-{R}aphson consensus for distributed convex optimization, {IEEE}
  Trans. on Automatic Control 61~(4) (2016) 994--1009.

\bibitem{nedic2017achieving}
A.~Nedi{\'c}, A.~Olshevsky, W.~Shi, Achieving geometric convergence for
  distributed optimization over time-varying graphs, {SIAM} J. on Optimization
  27~(4) (2017) 2597--2633.

\bibitem{qu2018harnessing}
G.~Qu, N.~Li, Harnessing smoothness to accelerate distributed optimization,
  {IEEE} Trans. on Control of Network Systems 5~(3) (2018) 1245--1260.

\bibitem{xu2018convergence}
J.~Xu, S.~Zhu, Y.~C. Soh, L.~Xie, Convergence of asynchronous distributed
  gradient methods over stochastic networks, {IEEE} Trans. on Autom. Control
  63~(2) (2018) 434--448.

\bibitem{xi2018addopt}
C.~Xi, R.~Xin, U.~A. Khan, {ADD-OPT}: Accelerated distributed directed
  optimization, {IEEE} Trans. on Autom. Control 63~(5) (2018) 1329--1339.

\bibitem{necoara2015linear}
I.~Necoara, V.~Nedelcu, On linear convergence of a distributed dual gradient
  algorithm for linearly constrained separable convex problems, Automatica 55
  (2015) 209--216.

\bibitem{alghunaim2018dual}
S.~Alghunaim, K.~Yuan, A.~Sayed, Dual coupled diffusion for distributed
  optimization with affine constraints, in: IEEE Conference on Decision and
  Control (CDC), 2018, pp. 829--834.

\bibitem{chang2014distributed}
T.-H. Chang, A.~Nedi{\'c}, A.~Scaglione, Distributed constrained optimization
  by consensus-based primal-dual perturbation method, {IEEE} Trans. on
  Automatic Control 59~(6) (2014) 1524--1538.

\bibitem{mateos2017distributed}
D.~Mateos-N{\'u}nez, J.~Cort{\'e}s, Distributed saddle-point subgradient
  algorithms with {L}aplacian averaging, {IEEE} Transactions on Automatic
  Control 62~(6) (2017) 2720--2735.

\bibitem{simonetto2016primal}
A.~Simonetto, H.~Jamali-Rad, Primal recovery from consensus-based dual
  decomposition for distributed convex optimization, Journal of Optimization
  Theory and Applications 168~(1) (2016) 172--197.

\bibitem{falsone2017dual}
A.~Falsone, K.~Margellos, S.~Garatti, M.~Prandini, Dual decomposition for
  multi-agent distributed optimization with coupling constraints, Automatica 84
  (2017) 149--158.

\bibitem{chang2016proximal}
T.-H. Chang, A proximal dual consensus {ADMM} method for multi-agent
  constrained optimization, IEEE Trans. on Signal Processing 64~(14) (2016)
  3719--3734.

\bibitem{notarnicola2017constraint}
I.~Notarnicola, G.~Notarstefano, Constraint coupled distributed optimization: a
  relaxation and duality approach, preprint arXiv:1711.09221.

\bibitem{kia2017distributed}
S.~S. Kia, Distributed optimal in-network resource allocation algorithm design
  via a control theoretic approach, Systems \& Control Letters 107 (2017)
  49--57.

\bibitem{zhang2018consensus}
Y.~Zhang, M.~M. Zavlanos, A consensus-based distributed augmented {L}agrangian
  method, in: IEEE Conference on Decision and Control (CDC), 2018, pp.
  1763--1768.

\bibitem{vujanic2016decomposition}
R.~Vujanic, P.~M. Esfahani, P.~J. Goulart, S.~Mari{\'e}thoz, M.~Morari, A
  decomposition method for large scale {MILP}s, with performance guarantees and
  a power system application, Automatica 67 (2016) 144--156.

\end{thebibliography}

\appendix
\renewcommand{\theequation}{A.\arabic{equation}}

\section{Proofs} \label{sec:appendix}

\subsection*{Proof of Lemma~\ref{lem:tracking_property}}
We prove \eqref{eq:avg_tracker_prop} by induction. For $\iter = 0$, given the initialization in
Step~\ref{step:tracker_init}, we have that
\begin{align*}
    \bar{d}_0 = \frac{1}{\Nag} \sum_{i = 1}^\Nag d_{i,0}
        & = \frac{1}{\Nag} \sum_{i = 1}^\Nag \left( A_i x_{i,0} - b_i \right) \\
        & = \frac{1}{\Nag} \left ( \sum_{i = 1}^\Nag A_i x_{i,0} - b \right ) 
\end{align*}
which proves that \eqref{eq:avg_tracker_prop} holds for $\iter = 0$. Assume now that \eqref{eq:avg_tracker_prop} holds up to $\iter$. If we can prove that \eqref{eq:avg_tracker_prop} holds also for $\iter + 1$, then the proof is completed. Indeed,
\begin{align*}
    \dkpavg &= \frac{1}{\Nag} \sum_{i = 1}^\Nag \dikp \\
        & \eqlabel{(a)} \frac{1}{\Nag} \sum_{i = 1}^\Nag \Big( \sum_{j = 1}^\Nag w_{ij} \djk + A_i \xikp - A_i \xik \Big) \\
        & \eqlabel{(b)} \frac{1}{\Nag} \sum_{j = 1}^\Nag \Big( \sum_{i = 1}^\Nag w_{ij} \Big) \djk
        + \frac{1}{\Nag} \sum_{i = 1}^\Nag \left( A_i \xikp - A_i \xik \right) \\
        & \eqlabel{(c)} \dkavg + \frac{1}{\Nag} \sum_{i = 1}^\Nag \left( A_i \xikp - A_i \xik \right) \\
        & \eqlabel{(d)} \frac{1}{\Nag} \Big( \sum_{i=1}^\Nag A_i \xikp - b \Big),
\end{align*}
where in (a) we used Steps~\ref{step:tracker_update} and~\ref{step:tracker_consensus} in (b) we exchanged the summations and we used column stochasticity of the weights (cf. Assumption~\ref{ass:balance_info}), in (c) we used the definition of $\dkavg$ and the induction step in (d). \qed

\subsection*{Proof of Lemma~\ref{lem:avg_dual_prop}}
By definition of $\lkpavg$, we can write
\begin{equation*}
    \lkpavg = \frac{1}{\Nag} \sum_{i = 1}^\Nag \likp
\end{equation*}
which can be further elaborated as follows
\begin{align*}
    \lkpavg & \eqlabel{(a)} \frac{1}{\Nag} \sum_{i = 1}^\Nag \bigg( \sum_{j = 1}^\Nag w_{ij} \ljk + c \, \dikp \bigg) \\
        & \eqlabel{(b)} \frac{1}{\Nag} \sum_{j = 1}^\Nag \left ( \sum_{i = 1}^\Nag w_{ij} \right ) \ljk + c \frac{1}{\Nag} \sum_{i = 1}^\Nag \dikp \\
        & \eqlabel{(c)} \lkavg + c \, \dkpavg,
\end{align*}
where in (a) we used the definition of $\likp$ in Step~\ref{step:dual_update}, in (b) we switched the summations, and in (c) we used column stochasticity of the weights (cf. Assumption~\ref{ass:balance_info}) together with the definition of $\lkavg$ and $\dkpavg$. \qed

\subsection*{Proof of Lemma~\ref{lem:bounded_errors}}
By \eqref{eq:all_primal_update} we have that $\xk \in X$ for all $\iter$ and, under Assumption~\ref{ass:convexity}, $X$ is bounded. Therefore the sequence $\seq{\xk}$ is bounded. Under Assumption~\ref{ass:balance_info}, from \eqref{eq:avg_tracker_prop} in Lemma~\ref{lem:tracking_property} we have that $\dkavg = (1/\Nag)(A\xk - b)$ and, therefore, also $\seq{\ddkavg}$ is bounded since $\ddkavg = \bone_\Nag \kron \dkavg$. Moreover, since $\zk = A_d \xk - \ddkavg$, then also $\seq{\zk}$ is bounded. This proves point~\ref{item:x_z_bounded}.

Being the sequence $\seq{\zk}$ bounded, we have that $\seq{\zkp - \zk}$ is also bounded. Therefore, exploiting the input-to-state stability property of the dynamical system in~\eqref{eq:as_stable_ed_dyn} under Assumptions~\ref{ass:connectivity} and~\ref{ass:balance_info}, we can conclude that $\seq{\edk}$ is a bounded sequence. Finally, owing to the input-to-state stability of the dynamical system in~\eqref{eq:as_stable_el_dyn} under Assumptions~\ref{ass:connectivity} and~\ref{ass:balance_info} together with boundedness of $\seq{\edk}$, we also have that $\seq{\elk}$ is bounded, thus proving point~\ref{item:ed_el_bounded} and concluding the proof. \qed

\subsection*{Proof of Proposition~\ref{prop:optimality_descent}}
By \cite[Lemma~4.1]{bertsekas1989parallel} applied to \eqref{eq:all_primal_update}, we have that $\xkp$ satisfies
\begin{multline*}
    f(\xkp) + [W\lk + c(A_d \xkp - A_d \xk + W\dk)]\T A_d \xkp \\
        \leq f(\xx) + [W\lk + c(A_d \xkp - A_d \xk + W\dk)]\T A_d \xx,
\end{multline*}
for any $\xx \in X$.
Recalling that $\dkp = W\dk + A_d \xkp - A_d \xk$ and $\lkp = W\lk + c\,\dkp$ (cf. \eqref{eq:all_tracker_update} and~\eqref{eq:all_dual_update}, respectively), we have
\begin{align*}
    f(\xkp) + \lkp\T A_d \xkp \leq f(\xx) + \lkp\T A_d \xx,
\end{align*}
for any $\xx \in X$. Under Assumption~\ref{ass:saddle}, we can set $\xx$ equal to some optimal solution $\xs$ of~\ref{eq:problem}, to obtain
\begin{align*}
    f(\xkp) + \lkp\T [A_d \xkp - A_d \xs] \leq f(\xs).
\end{align*}
Under Assumption~\ref{ass:saddle}, by the Saddle Point Theorem in \cite[pag. 665]{bertsekas1989parallel} we have that
\begin{equation} \label{eq:saddle_ineq}
    L(\xs,\lambda) \leq L(\xs,\ls) \leq L(\xx,\ls),
\end{equation}
for all $\xx \in X$ and for all $\lambda$. Given~\eqref{eq:saddle_ineq} and the fact that $L(\xs,\ls) = f(\xs)$, we have that
\begin{align*}
    f(\xkp) + \lkp\T [A_d \xkp - A_d \xs] &\leq f(\xs) \\
        &\leq f(\xx) + {\ls}\T [A \xx - b],
\end{align*}
for all $\xx \in X$ and any $\ls$. Writing $b = A \xs$ and noticing that ${\ls}\T [A \xx - A \xs] = {\lls}\T [A_d \xx - A_d \xs]$, it holds
\begin{multline*}
    f(\xkp) + \lkp\T [A_d \xkp - A_d \xs] \\
        \leq f(\xx) + {\lls}\T [A_d \xx - A_d \xs].
\end{multline*}
Setting $\xx = \xkp$, simplifying $f(\xkp)$, and bringing everything on the left-hand side, we have
\begin{align*}
    [\lkp - \lls]\T [A_d \xkp - A_d \xs] \leq 0.
\end{align*}
Adding and subtracting $\llkpavg\T [A_d \xkp - A_d \xs]$ yields
\begin{equation} \label{eq:opt_inequality1}
    \begin{multlined}
    [\llkpavg - \lls]\T [A_d \xkp - A_d \xs] \\
        + [\lkp - \llkpavg]\T [A_d \xkp - A_d \xs] \leq 0.
    \end{multlined}
\end{equation}
Focusing on the first term in \eqref{eq:opt_inequality1}, we have
\begin{align}
    [\llkpavg - \lls]\T [A_d &\xkp - A_d \xs] \nonumber \\
        &\eqlabel{(a)} \sum_{i = 1}^\Nag [\lkpavg - \ls]\T [A_i \xikp - A_i \xis] \nonumber \\
        &\eqlabel{(b)} [\lkpavg - \ls]\T \sum_{i = 1}^\Nag [A_i \xikp - A_i \xis] \nonumber \\
        &\eqlabel{(c)} [\lkpavg - \ls]\T [A \xkp - A \xs] \nonumber \\
        &\eqlabel{(d)} [\lkpavg - \ls]\T [A \xkp - b] \nonumber \\
        &\eqlabel{(e)} [\lkpavg - \ls]\T \Nag \dkpavg \nonumber \\
        &\eqlabel{(f)} \frac{\Nag}{c} [\lkpavg - \ls]\T [\lkpavg - \lkavg] \nonumber \\
        &\eqlabel{(g)} \frac{1}{c} [\llkpavg - \lls]\T [\llkpavg - \llkavg], \label{eq:opt_term1}
\end{align}
where (a) is due to the block structure of the vectors involved in the product, (b) is given by the fact that $[\lkpavg - \ls]$ does not depend on $i$, (c) is given by the definition of $A$, (d) is given by the fact that $A \xs = b$, (e) is due to \eqref{eq:avg_tracker_prop} in Lemma~\ref{lem:tracking_property}, (f) is given by \eqref{eq:avg_dual_prop} in Lemma~\ref{lem:avg_dual_prop}, and (g) uses the definitions of $\llkavg$ and $\lls$.
Owing to the fact that
\begin{align*}
    [\lkp - \llkpavg]\T \ddkpavg &= \sum_{i = 1}^{\Nag} [\likp - \lkpavg]\T \dkpavg \\
        &= \bigg( \sum_{i = 1}^{\Nag} [\likp - \lkpavg] \bigg)\T \dkpavg \\
        &= [\Nag \lkpavg - \Nag\lkpavg]\T \dkpavg \\
        &= 0
\end{align*}
and considering the second term in \eqref{eq:opt_inequality1}, we have
\begin{align}
    &[\lkp - \llkpavg]\T [A_d \xkp - A_d \xs] \nonumber \\
        &\quad= [\lkp - \llkpavg]\T [A_d \xkp - \ddkpavg - A_d \xs] \nonumber \\
        &\quad \eqlabel{(a)} [\lkp - \llkpavg]\T [\zkp - A_d \xs] \label{eq:opt_term2}
\end{align}
where in (a) we used the definition of $\zkp$ (cf.~\eqref{eq:z_def}). Using \eqref{eq:opt_term1} and \eqref{eq:opt_term2} in \eqref{eq:opt_inequality1} and multiplying by $2c$ we get
\begin{equation} \label{eq:opt_inequality2}
    \begin{multlined}
        2 [\llkpavg - \lls]\T [\llkpavg - \llkavg] \\
            + 2 c [\lkp - \llkpavg]\T [\zkp - A_d \xs] \leq 0.
    \end{multlined}
\end{equation}
Finally, using
\begin{equation*}
    2[a-b]\T [a-c] = \norm{a-b}^2 + \norm{a-c}^2 - \norm{c-b}^2
\end{equation*}
with the identifications $a = \llkpavg$, $b = \lls$ and $c = \llkavg$, \eqref{eq:optimality_descent} follows, thus concluding the proof. \qed

\subsection*{Proof of Theorem~\ref{thm:convergence}}
The dynamical evolution of $\edk$ and $\elk$, given in equations~\eqref{eq:as_stable_ed_dyn}
and~\eqref{eq:as_stable_el_dyn}, respectively, can be arranged as
\begin{equation}
    \ekp = F \ek + G(\ukp - \uk), \label{eq:consensus_dynamics}
\end{equation}
with
\begin{equation*}
    \ek = \begin{bmatrix} \elk \\ c\,\edk \end{bmatrix}\!,
    \quad
    F = \begin{bmatrix} \tW & \tW \\ 0 & \tW \end{bmatrix}\!,
    \quad
    G = \begin{bmatrix} I \\ I \end{bmatrix}\!,
\end{equation*}
and $\uk = c [ \zk - A_d \xs ]$. Let $P = P\T \succ 0$ and consider
\begin{align}
    \norm{G\ukp - \ekp}_{P}^2
    & \eqlabel{(a)} \norm{G\uk - \ek + (I-F) \ek}_{P}^2 \nonumber \\
        & \eqlabel{(b)} \norm{G\uk - \ek}_{P}^2 \nonumber \\
            &\qquad+ 2 [G\uk - \ek]\T P (I-F) \ek \nonumber \\
            &\qquad+ \ek\T (I-F\T) P (I-F) \ek \nonumber \\
        & \eqlabel{(c)} \norm{G\uk - \ek}_{P}^2 - \norm{\ek}_{P - F\T P F}^2 \nonumber \\
            &\qquad+ 2 \uk\T G\T P (I-F) \ek, \label{eq:lyapunov_consensus}
\end{align}
where in (a) we used the dynamics \eqref{eq:consensus_dynamics} and added $\pm \ek$ inside the norm, in (b) we expanded the square, while in (c) we performed some algebraic manipulations.
Summing~\eqref{eq:lyapunov_consensus} to the optimality-based inequality~\eqref{eq:optimality_descent} (cf. Proposition~\ref{prop:optimality_descent}) yields
\begin{align}
    & \norm{\llkpavg - \lls}^2 + 2 \ukp\T H \ekp + \norm{G\ukp - \ekp}_{P}^2 \nonumber \\
        & \quad\leq \norm{\llkavg - \lls}^2 + 2 \uk\T G\T P (I-F) \ek + \norm{G\uk - \ek}_{P}^2 \nonumber \\
            & \quad\qquad- \norm{\llkpavg - \llkavg}^2 - \norm{\ek}_{P - F\T P F}^2, \label{eq:lyapunov_combined}
\end{align}
with $H = [I \quad 0]$. Notice that, choosing $P$ such that $P - F\T P F \succ 0$, if the terms at $\iter+1$ were analogous to the terms at $\iter$, then~\eqref{eq:lyapunov_combined} would describe a descent condition. To enforce this, we need to find a $P$ satisfying the following conditions
\begin{subequations} \label{eq:lyapunov_conditions}
    \begin{align}
        &H = G\T P (I-F), \label{eq:cross_term} \\
        &P - F\T P F \succ 0, \label{eq:Q_psd} \\
        &P \succ 0. \label{eq:P_psd}
    \end{align}
\end{subequations}
Satisfying \eqref{eq:lyapunov_conditions} translates into~\eqref{eq:lyapunov_combined} being a non-expansive condition for the sequence $\seq{\norm{\llkpavg - \lls}^2 + 2 \ukp\T H \ekp + \norm{G\ukp - \ekp}_{P}^2}$, i.e., a stability condition for the nonlinear dynamical system representing Algorithm~\ref{alg:DT-ADMM}.

Take $P$ partitioned in blocks as follows
\begin{equation*}
    P = \begin{bmatrix} P_1 & P_2 \\ P_2\T & P_3 \end{bmatrix}.
\end{equation*}
Owing to the fact that all eigenvalues of $\tW$ lies in the open unit circle, thus $(I - \tW)$ is invertible and, from the equality $M M^{-1} = M^{-1} M = I$ with $M = (I - \tW)$, we have the following identities
\begin{subequations} \label{eq:tW_identities}
    \begin{align}
    & \tW (I - \tW)^{-1} = (I - \tW)^{-1} - I, \label{eq:tW_lid} \\
    & (I - \tW)^{-1} \tW = (I - \tW)^{-1} - I. \label{eq:tW_rid}
    \end{align}
\end{subequations}
Noticing that
\begin{align*}
    (I - F)^{-1} = \begin{bmatrix} (I - \tW)^{-1} & (I - \tW)^{-1} \tW (I - \tW)^{-1} \\ 0 & (I - \tW)^{-1}\end{bmatrix},
\end{align*}
condition \eqref{eq:cross_term} can be rewritten as $G\T P = H (I - F)^{-1}$ and translates into the following constraints on the blocks of $P$:
\begin{subequations} \label{eq:P_relations}
    \begin{align}
        P_1 + P_2\T &= (I - \tW)^{-1} \label{eq:P12_rel} \\
        P_2 + P_3 &= (I - \tW)^{-1} \tW (I - \tW)^{-1} \nonumber \\
            & \eqlabel{(a)} (I - \tW)^{-2} - (I - \tW)^{-1} \label{eq:P23_rel} \\
        P_1 + P_2\T + P_2 + P_3 &= (I - \tW)^{-2}, \label{eq:P123_rel}
    \end{align}
\end{subequations}
where in (a) we used~\eqref{eq:tW_lid}, while~\eqref{eq:P123_rel} follows by simply summing~\eqref{eq:P12_rel} and~\eqref{eq:P23_rel} and will be useful in the forthcoming discussion. Note that, from \eqref{eq:P12_rel}, we can further deduce that $P_2 = P_2\T$ since both $P_1$ and $(I - \tW)^{-1}$ are symmetric matrices owing to $P$ being symmetric and Assumption~\ref{ass:balance_info}, respectively.

We can thus express the matrix $P - F\T P F$ as
\begin{align}
    &
    \begin{bmatrix}
        P_1 - \tW P_1 \tW                   &P_2 - \tW (P_1 + P_2) \tW \\
        P_2\T - \tW (P_1 + P_2\T) \tW       &P_3 - \tW (P_1 + P_2\T + P_2 + P_3) \tW
    \end{bmatrix}
    \nonumber \\
    &\quad \eqlabel{(a)}
    \begin{bmatrix}
        P_1 - \tW P_1 \tW                   &P_2 - \tW (I - \tW)^{-1} \tW \\
        P_2\T - \tW (I - \tW)^{-1} \tW      &P_3 - \tW (I - \tW)^{-2} \tW
    \end{bmatrix}
    \nonumber \\
    &\quad \eqlabel{(b)}
    \begin{bmatrix}
        P_1 - \tW P_1 \tW                   & \tW - (P_1 - I) \\
        \tW - (P_1 - I)                     &P_1 - I
    \end{bmatrix},
\end{align}
where in (a) we used~\eqref{eq:P_relations} and in (b) we leveraged identities in~\eqref{eq:tW_identities} together with \eqref{eq:P12_rel} and \eqref{eq:P23_rel} to express $P_2$ and $P_3$ as a function of $P_1$ only.

Next, we find a value for $P_1 \succ 0$ ensuring that $P - F\T P F$ is positive definite. This requirement can be posed in terms of its blocks by means of the Schur complement, i.e., $P - F\T P F \succ 0$ if and only if
\begin{subequations} \label{eq:Schur_Qpsd}
    \begin{align}
        & P_1 - I \succ 0, \label{eq:Schur_condition1_Qpsd} \\
        & P_1 - \tW P_1 \tW - (\tW - (P_1 - I)) (P_1 - I)^{-1} (\tW - (P_1 - I)) \nonumber \\
        & \hspace{1em}= 2\tW + I - \tW (P_1 + (P_1 - I)^{-1}) \tW \succ 0. \label{eq:Schur_condition2_Qpsd}
    \end{align}
\end{subequations}
Take $P_1 = 2I$ and let $T \Lambda T\T = \tW$ be the eigendecomposition of the symmetric matrix $\tW$, with $T T\T = I$ and $\Lambda$ being a diagonal matrix containing the eigenvalues of $\tW$. Condition \eqref{eq:Schur_condition1_Qpsd} is trivially satisfied by this choice of $P_1$, while \eqref{eq:Schur_condition2_Qpsd} becomes
\begin{align*}
    2\tW + I - \tW (2I + I) \tW &= 2\tW + I - 3 \tW^2 \\
        &= T (I + 2\Lambda - 3\Lambda^2) T\T \succ 0,
\end{align*}
which is equivalent to $I + 2\Lambda - 3\Lambda^2 \succ 0$ and always holds true when all the eigenvalues of $\tW$ lies within $(-1/3,1)$. In view of Assumption~\ref{ass:W_eigval}, this condition is guaranteed.

Finally, using again Schur's complement, it is also easy to prove that $P_1 = 2I$ satisfies condition~\eqref{eq:P_psd}. In fact, we have that
\begin{equation*}
    P = \begin{bmatrix} 2I & (I \!-\! \tW)^{-1} - 2I \\ (I \!-\! \tW)^{-1} - 2I & (I \!-\! \tW)^{-2} -2(I \!-\! \tW)^{-1} + 2I \end{bmatrix} \succ 0
\end{equation*}
if and only if
\begin{subequations} \label{eq:Schur_Ppsd}
    \begin{align}
        &2I \succ 0, \label{eq:Schur_condition1_Ppsd} \\
        &(I \!-\! \tW)^{-2} -2(I \!-\! \tW)^{-1} + 2I - \frac{1}{2}\left( (I \!-\! \tW)^{-1} - 2I \right)^2 \nonumber \\
        &\qquad= \frac{1}{2} (I - \tW)^{-2} \succ 0. \label{eq:Schur_condition2_Ppsd}
    \end{align}
\end{subequations}
Condition~\eqref{eq:Schur_condition2_Ppsd} is satisfied since $\tW$ has all eigenvalues in the open unit circle and hence $I - \tW \succ 0$.

From \eqref{eq:cross_term} we can rewrite \eqref{eq:lyapunov_combined} as
\begin{align}
    &\norm{\llkpavg - \lls}^2 + 2 \ukp\T H \ekp + \norm{G\ukp - \ekp}_{P}^2 \nonumber \\
        &\quad\leq \norm{\llkavg - \lls}^2 + 2 \uk\T H \ek + \norm{G\uk - \ek}_{P}^2 \nonumber \\
            &\quad\qquad- \norm{\llkpavg - \llkavg}^2 - \norm{\ek}_{Q}^2, \label{eq:lyapunov_final}
\end{align}
where we set $Q = P - F\T P F \succ 0$.
Rearranging some terms, and summing \eqref{eq:lyapunov_final} from $\iter = 0$ to $M-1$, with $M\in\bN$, we have
\begin{align*}
    \sum_{\iter = 0}^{M-1} &\norm{\llkpavg - \llkavg}^2 + \norm{\ek}_{Q}^2 \\
        &\qquad \leq \norm{\bs{\bar{\ll}}_0 - \lls}^2 + 2 \uu_0\T H \ee_0 + \norm{G \uu_0 - \ee_0}_{P}^2 \\
        &\qquad\quad - \norm{\bs{\bar{\ll}}_M - \lls}^2 - 2 \uu_M\T H \ee_M - \norm{G \uu_M - \ee_M}_{P}^2 \\
        &\qquad \leq \norm{\bs{\bar{\ll}}_0 - \lls}^2 + \norm{G \uu_0 - \ee_0}_{P}^2 + 2\cC,
\end{align*}
where in the second inequality we neglected the terms $-\norm{\bs{\bar{\ll}}_M - \lls}^2$ and $-\norm{G \uu_M - \ee_M}_{P}^2$ in the right hand side since they are non-positive (recall $P \succ 0$) and we used the fact that $|2 \uk\T H \ek| \leq \cC$ for all $\iter \geq 0$ owing to the results of Lemma~\ref{lem:bounded_errors}. Taking the limit as $M \to \infty$ finally yields
\begin{equation*}
    \sum_{\iter = 0}^\infty \left ( \norm{\llkpavg - \llkavg}^2 + \norm{\ek}_{Q}^2 \right ) < \infty,
\end{equation*}
which, recalling that $Q \succ 0$, implies that $\lim_{\iter \to \infty} \norm{\llkpavg - \llkavg} = 0$ and $\lim_{\iter \to \infty} \norm{\ek} = 0$. Hence, $\lim_{\iter \to \infty} \norm{\edk} = 0$ and $\lim_{\iter \to \infty} \norm{\elk} = 0$, thus proving points~\ref{item:thm_tracker_consensus} and~\ref{item:thm_dual_consensus}, respectively.

Since $\ddkpavg = \frac{1}{c} (\llkpavg - \llkavg)$, we also have that $\lim_{\iter \to \infty} \norm{\ddkavg} = 0$ and therefore $\lim_{\iter \to \infty} \norm{\dkavg} = 0$, thus proving point~\ref{item:thm_feasibility}.

From \eqref{eq:lyapunov_final} we also have that the sequence
\begin{equation*}
    \seq{ \norm{\llkavg - \lls}^2 + 2 \uk\T H \ek + \norm{G\uk - \ek}_{P}^2 }
\end{equation*}
is non-increasing, bounded below since $|2 \uk\T H \ek| \leq \cC$ and $\norm{G\uk - \ek}_{P}^2 \geq 0$ (recall $P \succ 0$), and, therefore, convergent. Since $\norm{\ek}$ is vanishing and $\seq{\uk}$ is bounded, we have that also the sequence $\seq{ \norm{\llkavg - \lls}^2 + \norm{G\uk}_{P}^2 }$ is convergent.
A straightforward computation using \eqref{eq:P123_rel} shows that $G\T P G = (I - \tW)^{-2} \succ 0$ which, recalling the definition of $\uk = c\,(\zk - A_d \xs)$, implies that the sequence
\begin{equation*}
    \seq{ \norm{\llkavg - \lls}^2 + c^2 \norm{\zk - A_d \xs}_{(I - \tW)^{-2}}^2 }
\end{equation*}
is convergent, which finally implies that
\begin{equation*}
    \seq{ \norm{\llkavg - \lls}^2 + c^2 \norm{\zk - A_d \xs}^2 }
\end{equation*}
is convergent due to norm equivalence, thus proving point~\ref{item:thm_convergence} and concluding the proof. \qed

\subsection*{Proof of Theorem~\ref{thm:main_result}}
Following the first steps in the proof of Proposition~\ref{prop:optimality_descent}, for any optimal solution $\xs$ for~\ref{eq:problem}, we can write
\begin{align*}
    f(\xkp) + \lkp\T [A_d \xkp - A_d \xs] \leq f(\xs).
\end{align*}
By adding and subtracting $\llkpavg\T [A_d \xkp - A_d \xs]$ to the left-hand side, we have
\begin{multline*}
    f(\xkp) + \llkpavg\T [A_d \xkp - A_d \xs] \\
        + {\elkp}\T [A_d \xkp - A_d \xs] \leq f(\xs).
\end{multline*}
Noticing that
\begin{align}
    \llkpavg\T [A_d \xkp - A_d \xs] &= \lkpavg\T [A \xkp - A \xs] \nonumber \\
        &= \lkpavg\T [A \xkp - b] \nonumber \\
        &= \lkpavg\T \, \Nag \,\dkpavg, \label{eq:lldd_ld_equality}
\end{align}
where the last inequality is due to Lemma~\ref{lem:tracking_property}, we have that
\begin{equation*}
    f(\xkp) + \Nag \lkpavg\T \dkpavg + {\elkp}\T [A_d \xkp - A_d \xs] \leq f(\xs).
\end{equation*}
From Theorem~\ref{thm:convergence}\ref{item:thm_convergence} we have that the sequence
\begin{equation*}
    \seq{ \norm{\llkavg - \lls}^2 + c^2 \norm{\zk - A_d \xs}^2 }
\end{equation*}
is bounded. Since $\seq{\zk}$ is also bounded by Lemma~\ref{lem:bounded_errors}, we have that $\seq{\llkavg}$ and, therefore, also $\seq{\lkavg}$ are both bounded. Recalling that $\seq{\xk}$ is bounded, by Theorem~\ref{thm:convergence} points~\ref{item:thm_feasibility} and~\ref{item:thm_dual_consensus} respectively, it holds
\begin{align}
    &\lim_{\iter \to \infty} \lkpavg\T \dkpavg = 0, \label{eq:limit_cross1} \\
    &\lim_{\iter \to \infty} {\elkp}\T [A_d \xkp - A_d \xs] = 0, \label{eq:limit_cross2}
\end{align}
and therefore
\begin{align}
    \limsup_{\iter \to \infty} f&(\xkp) + \Nag \lkpavg\T \dkpavg + {\elkp}\T [A_d \xkp - A_d \xs] \nonumber \\
        &= \limsup_{\iter \to \infty} f(\xkp) \leq f(\xs). \label{eq:limsup_xk}
\end{align}
Since $\xk \in X$ for all $\iter \geq 0$ and, by Lemma~\ref{lem:tracking_property} and Theorem~\ref{thm:convergence}\ref{item:thm_feasibility},
\begin{equation*}
    \lim_{\iter \to \infty} \norm{A\xk -b} = \lim_{\iter \to \infty} \norm{\dkavg} = 0,
\end{equation*}
we know that any limit point $\xlim$ of the sequence $\seq{\xk}$ is feasible for~\ref{eq:problem} and, thus, satisfies $f(\xlim) \geq f(\xs)$. Combining this fact with \eqref{eq:limsup_xk} we can conclude that $f(\xlim) = f(\xs)$ for any limit point $\xlim$. Since all limit points of $\seq{\xk}$ are feasible and achieve the optimal value, they are optimal solutions for~\ref{eq:problem}, proving
Theorem~\ref{thm:main_result}\ref{eq:thm_primal_opt}.

Next, we prove the convergence of the local dual variables estimates to an optimal dual solution. Similarly to the proof of Proposition~\ref{prop:optimality_descent}, by applying \cite[Lemma~4.1]{bertsekas1989parallel} to Step~\ref{step:primal_update} we have that $\xikp$ satisfies
\begin{multline*}
    f_i(\xikp) + [\ellik + c(A_i \xikp - A_i \xik + \deltaik)]\T A_i \xikp \\
        \leq f_i(x_i) + [\ellik + c(A_i \xikp - A_i \xik + \deltaik)]\T A_i x_i,
\end{multline*}
for all $x_i \in X_i$. Noticing that $\dikp = \deltaik + A_i \xikp - A_i \xik$ and $\likp = \ellik + c\,\dikp$ from Steps~\ref{step:tracker_update} and~\ref{step:dual_update} respectively and subtracting $\likp\T b_i$, we have
\begin{multline*}
    f_i(\xikp) + \likp\T (A_i \xikp - b_i) \\
        \leq f_i(x_i) + \likp\T (A_i x_i - b_i)
\end{multline*}
for all $x_i \in X_i$, which is equivalent to
\begin{equation*}
    \xikp \in \argmin_{x_i \in X_i} f_i (x_i) + \likp\T ( A_i x_i - b_i)
\end{equation*}
and, thus, by definition of $\dual_i(\lambda)$ in \eqref{eq:local_dual_function}, we have
\begin{equation*}
    \dual_i(\likp) = f_i(\xikp) + \likp\T (A_i \xikp - b_i),
\end{equation*}
if $b_1,\dots,b_\Nag$ are such that $\sum_{i=1}^\Nag b_i = b$. Setting $b_i = A_i \xis$, summing over $i = \toN$, adding and subtracting $\llkpavg\T [A_d \xkp - A_d \xs]$ to the right hand side, and exploiting again \eqref{eq:lldd_ld_equality}, we obtain
\begin{equation}\label{eq:dual_cost_iter}
    \begin{multlined}
        \sum_{i=1}^\Nag \dual_i(\likp) = f(\xkp) + \Nag \, \lkpavg\T \dkpavg \\
            + {\elkp}\T [A_d \xkp - A_d \xs].
    \end{multlined}
\end{equation}
Since the sequence $\seq{\lkpavg}$ is bounded, it admits a convergent subsequence $\subseq{(\lkpavg,\xkp)}$ with $\cK \subseteq \bN$. Let $(\lavglim,\xlim)$ be its limit point. Now, taking the limit of \eqref{eq:dual_cost_iter} across $\cK$, we have
\begin{align}
    \lim_{\cK\ni \iter \to \infty} \sum_{i=1}^\Nag \dual_i(\likp)
        &= \lim_{\iter \in \cK} \Big\{ f(\xkp) + \Nag \, \lkpavg\T \dkpavg \nonumber \\
            &\qquad\qquad + {\elkp}\T [A_d \xkp - A_d \xs] \Big\} \nonumber \\
        & \eqlabel{(a)} f(\xlim) \nonumber \\
        & \eqlabel{(b)} f(\xs) \nonumber \\
        & \geqlabel{(c)} \max_{\lambda} \sum_{i=1}^\Nag \dual_i(\lambda), \label{eq:liminf_dual}
\end{align}
where in (a) we used \eqref{eq:limit_cross1} and \eqref{eq:limit_cross2}, the fact that limit points of $\seq{\xk}$ achieve the optimal value in (b), and weak duality in (c). On the other hand, from \eqref{eq:local_dual_function}, we have
\begin{align}
    \lim_{\cK\ni \iter \to \infty} \dual_i(\likp)
     & \leq \lim_{\cK\ni \iter \to \infty} \Big \{ f_i(x_i) + \likp\T (A_i x_i - b_i) \Big \} \nonumber \\
        &= f_i(x_i) + \Big ( \lim_{\cK\ni \iter \to \infty} \lkpavg \Big ) \T (A_i x_i - b_i) \nonumber \\
        &= f_i(x_i) + \lavglim\T (A_i x_i - b_i), \label{eq:limsup_dual1}
\end{align}
for all $x_i \in X_i$, where the first equality is due to $\lim_{\iter \to \infty} \lik = \lkavg$, $i = \toN $ (recall $\lim_{\iter \to \infty} \norm{\elk} = 0$), and in the last equality we plugged the limit point $\lavglim$ of $\subseq{\lkpavg}$. Setting $x_i$ in \eqref{eq:limsup_dual1} equal to the minimizer of $\dual_i(\lavglim)$ and summing over $i = \toN $ yields
\begin{equation*}
    \lim_{\cK\ni \iter \to \infty} \sum_{i=1}^\Nag \dual_i(\likp)
        \leq \sum_{i=1}^\Nag \dual_i(\lavglim) \leq \max_{\lambda} \sum_{i=1}^\Nag \dual_i(\lambda).
\end{equation*}
Combining the latter with \eqref{eq:liminf_dual} gives
\begin{equation} \label{eq:dual_optimality}
    \sum_{i=1}^\Nag \dual_i(\lavglim) = \max_{\lambda} \sum_{i=1}^\Nag \dual_i(\lambda),
\end{equation}
which means that $\lavglim$ is optimal for~\ref{eq:dual_problem}.

Recall that, from Theorem~\ref{thm:convergence}\ref{item:thm_convergence}, we have that the sequence
\begin{equation*}
    \seq{ \norm{\llkavg - \lls}^2 + c^2 \norm{\zk - A_d \xs}^2 }
\end{equation*}
is convergent for any pair $\xs$ and $\ls$ which are optimal for~\ref{eq:problem} and~\ref{eq:dual_problem}, respectively. We can thus take $\lls = \bone_\Nag \kron \lavglim$ and $\xs = \xlim$ and conclude that
\begin{equation*}
    \lim_{\cK\ni \iter \to \infty} \Big\{ \norm{\llkavg - \lls}^2 + c^2 \norm{\zk - A_d \xs}^2 \Big\} = 0.
\end{equation*}
But since $\seq{ \norm{\llkavg - \lls}^2 + c^2 \norm{\zk - A_d \xs}^2 }$ is convergent, all its limit points are the same, and therefore
\begin{equation*}
    \lim_{\cK\ni \iter \to \infty} \Big\{ \norm{\llkavg - \lls}^2 + c^2 \norm{\zk - A_d \xs}^2 \Big\} = 0,
\end{equation*}
for some primal optimal solution $\xs$ and some optimal dual solution $\ls$, which implies that, for all $i = \toN$,
\begin{align*}
    &\lim_{\iter \to \infty} \norm{\lik - \ls} = 0 \\
    &\lim_{\iter \to \infty} \norm{A_i \xik - A_i \xis} = 0
\end{align*}
since $\elk$ is vanishing. So that Theorem~\ref{thm:main_result}\ref{eq:thm_dual_opt} follows. \qed

\end{document}